\documentclass[12pt]{amsart}  
\usepackage{amsfonts,amssymb,graphicx}

\usepackage[normalem]{ulem} 

\usepackage{mathtools}\mathtoolsset{showonlyrefs,showmanualtags}

\usepackage{hyperref} 
\hypersetup{
    colorlinks=true,       
    linkcolor=blue,          
    citecolor=magenta,        
    filecolor=magenta,      
    urlcolor=cyan           
}

\usepackage{amsrefs}

\theoremstyle{plain} 
\newtheorem{lemma}[equation]{Lemma} 
\newtheorem{proposition}[equation]{Proposition} 
\newtheorem{theorem}[equation]{Theorem} 
\newtheorem{corollary}[equation]{Corollary}

\theoremstyle{definition}
\newtheorem{definition}[equation]{Definition} 

\theoremstyle{remark}
\newtheorem{remark}[equation]{Remark}

\numberwithin{equation}{section}

\begin{document}
\title[Weak and Strong $ A_p$ Estimates] {Weak And Strong Type Estimates for Maximal Truncations of Calder\'on-Zygmund Operators on $ A_p$ Weighted Spaces}

\author[Hyt\"onen]{Tuomas P.\ Hyt\"onen}
\address{Department of Mathematics and Statistics, University of Helsinki, Gustaf H\"allstr\"omin katu 2b, FI-00014 Helsinki, Finland}
\email{tuomas.hytonen@helsinki.fi}
\thanks{Research of TPH and HM supported by  the Academy of Finland grants 130166, 133264 and 218418}

\author[Lacey]{Michael T. Lacey}   
\thanks{Research of MTL, and MCR  supported in part by NSF grant 0968499}
\address{ School of Mathematics, Georgia Institute of Technology, Atlanta GA 30332, USA}
\email {lacey@math.gatech.edu}

\author[Martikainen]{Henri Martikainen}
\address{Department of Mathematics and Statistics, University of Helsinki, Gustaf H\"allstr\"omin katu 2b, FI-00014 Helsinki, Finland}
\email{henri.martikainen@helsinki.fi} 

\author[Orponen] {Tuomas Orponen}
\thanks{Research of TO supported by the Finnish Centre of Excellence in Analysis and Dynamics Research}
\address{Department of Mathematics and Statistics, University of Helsinki, Gustaf H\"allstr\"omin katu 2b, FI-00014 Helsinki, Finland}
\email{tuomas.orponen@helsinki.fi} 

\author[Reguera]{Maria  Carmen Reguera}
\address{ School of Mathematics, Georgia Institute of Technology, Atlanta GA 30332, USA}
\email {mreguera@math.gatech.edu}

\author[Sawyer]{Eric T. Sawyer}
\address{ Department of Mathematics \& Statistics, McMaster University, 1280
Main Street West, Hamilton, Ontario, Canada L8S 4K1 }
\email{sawyer@mcmaster.ca}
\thanks{Research of ETS supported in part by NSERC}

\author[Uriarte-Tuero]{Ignacio Uriarte-Tuero}
\address{ Department of Mathematics \\
Michigan State University \\
East Lansing MI }
\email{ignacio@math.msu.edu}
\thanks{Research of IU-T supported in part by the NSF, through grant DMS-0901524.}
\date{}

\begin{abstract}
For $ 1<p< \infty $, weight $ w\in A_p$ and any $ L ^2 $-bounded Calder\'on-Zygmund operator $ T$, we show that there is a constant $ C _{T,p}$ so that 
we have the weak and strong type inequalities 
\begin{align*}
\lVert T _{\natural} f \rVert_{ L ^{ p, \infty} (w) }& \le C _{T,p} \lVert w\rVert_{A_p} \lVert f\rVert_{ L ^{p} (w)}
\\
\lVert T _{\natural} f \rVert_{ L ^{ p} (w) } &\le C _{T,p} \lVert w\rVert_{A_p} ^{\max \{1, (p-1) ^{-1} \}} \lVert f\rVert_{ L ^{p} (w)} \,, 
\end{align*}
where $ T _{\natural} $ denotes the maximal truncations of $T$,  $ w$ is a weight, and $ \lVert w\rVert_{A_p}$ denotes the 
Muckenhoupt $ A_p$ characteristic of $ w$.  These estimates are not improvable in the power of $  \lVert w\rVert_{A_p}$. 
Our argument follows the outlines of the arguments of Lacey--Petermichl--Reguera (Math.\ Ann.\ 2010) and Hyt\"onen--P\'erez--Treil--Volberg (arXiv, 2010)
with new ingredients, including a weak-type estimate for certain duals of $ T _{\natural}$, and sufficient conditions for two weight inequalities in $ L ^{p}$ for $ T _{\natural}$. 
Our proof does not rely upon extrapolation.  
\end{abstract}

\maketitle

\setcounter{tocdepth}{1} 
\tableofcontents

\section{Overview and Introduction} 

Our subject is weighted inequalities for maximal truncations $ T _{\natural}$ of Calder\'on-Zygmund operators.  There are two main results.
 First, we prove weak and strong norm estimates on $ T _{\natural}$ on $ L ^{p} (w)$, that are  sharp in the $ A_p$ characteristic of the weight $ w$.   In the generality of this paper, this was only known for the untruncated operators, 
a question investigated by many, culminating in the definitive result in  \cite{1007.4330}.

Second, for dyadic Calder\'on-Zygmund operators, termed Haar Shift operators, we prove 
sufficient conditions for the weak and strong type two-weight inequalities $ T _{\natural}$. These estimates are effective in terms of 
a notion of complexity for the Haar shift, and while providing only sufficient conditions, are sharp enough in the $ A_p$ setting that we can 
conclude our first result from them.  

We recall definitions. 

\begin{definition}\label{d.czo} 
A \emph{Calder\'on-Zygmund  operator} in $\mathbb{R}^d$ is a \emph{bounded} in $L^2$ integral operator with kernel $K (x,y)$, 
defined by the expression 
\begin{equation*}
\langle T f ,g \rangle = \int\!\!\int 
f (x) g (y) K (x,y) \; dy 
\end{equation*}
for all continuous compactly supported functions $ f,g$ with $ \textup{dist} (\textup{supp} (f), \textup{supp} (g))>0$. 
The kernel $ K (x,y)$ satisfies  the following growth and  smoothness conditions for $ x, x', y, y'\in \mathbb R ^{d}$, $ x\neq y $
\begin{gather*}
|K(x, y)|   \le \frac{C _{T}}{|x-y|^d} \,, \qquad  x, y\in \mathbb{R}^d,\ x\ne y. 
\\
 |K(x, y) - K(x', y) | + |K(y, x) - K(y, x')| \le C_T \frac{|x-x'|^{\alpha}}{|x-y|^{d+\alpha}}\,, \qquad |x-x'|<|x-y|/2.
\end{gather*}
Here, $ C_T$ is an absolute constant.  We denote the maximal truncations of $ T$ by 
\begin{equation*}
T _{\natural} f (x) := \sup _{ 0<\epsilon <  \nu  } 
\Bigl\lvert \int _{\epsilon <  \lvert  y\rvert < \nu} f (y) K(x,y) \; dy    \Bigr\rvert
\end{equation*}
It is well-known that $ T$ and $ T _{\natural} $ extend to bounded operators on $ L ^{p} (\mathbb R ^{d})$, for $ 1<p< \infty $. 
\end{definition}

Prominent examples include the Hilbert and  Beurling transforms, as well as the vector $ R$ of Riesz transforms.  
If $ w$ is a weight on $ \mathbb R ^{d}$, namely a non-negative measure, with density also denoted as $ w$ that is non-negative almost everywhere, 
it is well-known \cite{MR0312139} that $ R$ is bounded on $ L ^{p} (w)$, $ 1<p< \infty $,  if and only if $ w$ satisfies the famous Muckenhoupt $ A_p$ condition 
\begin{equation}\label{e.ApDef}
\lVert w\rVert_{A_p} 
:= 
\sup _{Q}  \lvert  Q\rvert ^{-1} \int _{Q} w (dx) \Bigl[ \lvert  Q\rvert ^{-1}  \int _{Q} \sigma (dx) \Bigr] ^{p-1} 
\end{equation}
where $ \sigma $ is the weight with density $ w ^{-1/(p-1)}$, which is dual to $ w$.   Note that $ \lVert w\rVert_{A_p}$ is certainly not 
a norm. 

On the other hand, determining the sharp 
dependence of Calder\'on-Zygmund operators on the quantity $ \lVert w\rVert_{A_p}$ is not straight forward, as first pointed out 
by Buckley \cite{MR1124164}.   This direction has been intensively studied in recent years, with the sharp result for 
$ T$ established in \cite{1007.4330}, following the contributions of several.  We refer the reader to the introductions of 
\cites{1007.4330,1010.0755,MR2628851,MR2354322} for more information about the history and range of techniques 
brought to bear on this problem.  

Our first main result is this Theorem. 

\begin{theorem} \label{t.main}
\label{t.weakStrong}  For $  T $ an $ L ^{2} (\mathbb R ^{d})$ bounded  Calder\'on-Zygmund Operator, 
\begin{align}\label{e.weak}
\lVert   T _{\natural } f  \rVert _{ L ^{p,\infty } (w)} &\le C_T \lVert  w \rVert _{ A _{p}} \lVert  f \rVert _{L ^{p} (w)} \, 
\,,\qquad 1 < p < 2\,, 
\\ \label{e.strong}
\lVert   T _{\natural} f  \rVert _{ L ^{p } (w)} &\le C_T \lVert  w \rVert _{ A _{p}} ^{ \max \{1,  (p-1) ^{-1}\} } \lVert  f \rVert _{L ^{p} (w)} \, 
\,, \qquad 1 < p < \infty \,.
\end{align}
\end{theorem}

Well known examples involving power weights (see the conclusion of \cite{MR2354322}) show that all the estimates above 
are sharp. Indeed, these bounds match the best possible bounds for the untruncated operator $ T$.  
The weak-type estimate was conjectured by Andrei Lerner \cite{lerner}, who also conjectured that the maximal truncations 
should have the same behavior in the $ A_p$ charateristic as the untruncated operators (personal communication).  
As far as we are aware, this is the first place in which the sharp estimates for $ T _{\natural}$ have been established, and the weak-type inequality is new even for untruncated $T$. 

\begin{figure}
\begin{center}
\includegraphics{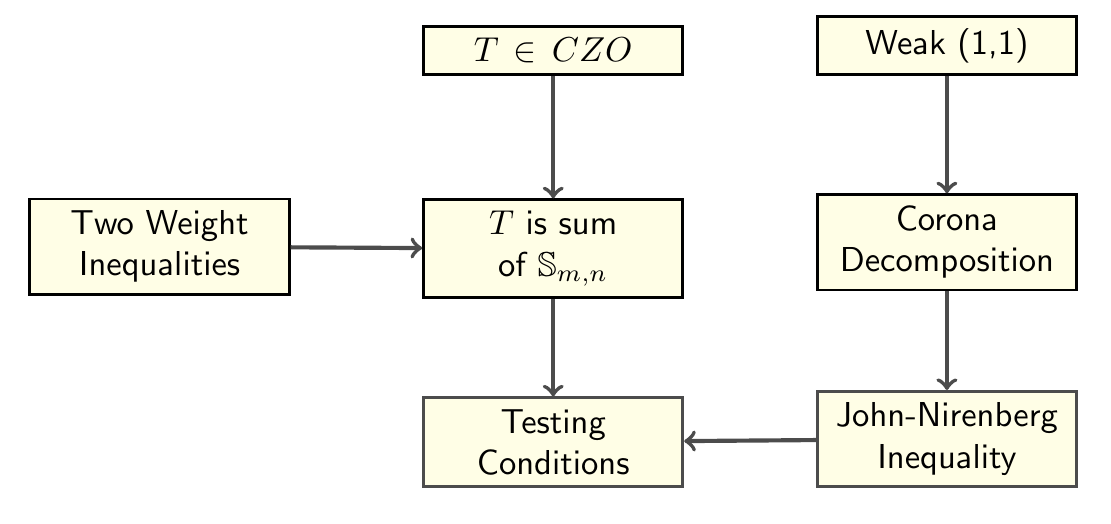} 
\end{center}
\caption{Proof Strategy}
\label{f.1}  
\end{figure}

We move to a discussion of the proof strategy for this Theorem.  We will follow the outlines of the argument of \cite{1010.0755}, but the underlying details are substantially different.  The strategy is summarized in Figure~\ref{f.1}, and has the following points.  

We begin with a Calder\'on-Zygmund Operator $ T$, and the important step, identified in \cite{1007.4330}, is to write $ T$ as a rapidly convergent sum of \emph{Haar Shift Operators} $ \mathbb S _{m,n}$.  See Definition~\ref{d.haarShift}, and Theorem~\ref{t.bcr}.  

 Haar Shift Operators are themselves dyadic variants of Calder\'on-Zygmund Operators, and come with an essential notion of 
\emph{complexity}, which is the measure of how many inter-related dyadic scales the operator has.   As Calder\'on-Zygmund Operators, they satisfy many estimates already, but it is a vital point that in order to use the fact that $ T$ is a rapidly convergent sum of these operators, all relevant estimates must be shown to be \emph{at most polynomial in complexity.}  We will refer to this as an \emph{effective estimate}. This requires that we revisit most facts 
about these operators, and verify that they meet this requirement. 

The next crucial stage, the most complicated part of this argument, is to prove reasonably sharp \emph{two weight inequalities} for Haar Shift Operators.  The import here is that to prove our theorem, 
 much of the argument must work in the generality of the two weight setting for a dyadic Calder\'on-Zygmund Operator.  
 That  the weight is in 
$ A_p$ is a fact that  can only be used very sparingly. In this, we are following  the pattern of \cites{MR2657437,1007.4330,0911.0713,1006.2530}.

 All of these prior works depended upon two-weight inequalities for the \emph{untruncated} operator, and \emph{only in $ L ^2 $.} 
Here, we are concerned with two-weight inequalities for the \emph{maximal truncations}; these estimates will apply in all $ L ^{p}$ spaces, an important point as concerns the weak-type inequality.  These estimates are taken up in \S\ref{s.2wt}, with the weak-type estimate being simple, and the strong type estimate being the most complicated estimate.   Different variants of this argument have been used in \cites{0911.3920,0807.0246,0911.0713,1006.2530}, with the point here being that the estimates in \S\ref{s.2wt} track complexity.  
See this section for more history on these estimates.

 The essential consequence of the two-weight inequalities is that they reduce the question of estimating the norm of $ T$ to that of testing the norm on a much simpler class of functions---weighted indicators of intervals.  These conditions are in turn verified by using a chain of arguments that begins with the verification of certain weak-$L^1$ inequalities for the Haar Shift Operators. We need these weak-type bounds for the \emph{adjoints of all linearizations} of the maximal truncation operator.
This is an estimate not of a classical nature, and is taken up in \S\ref{s.w11}.  

This weak-integrability has a certain measure of uniformity.  This permits the use of a \emph{John-Nirenberg Inequality} that shows that uniform weak-integrability actually implies exponential integrability.  This principle, again needed for certain maximal truncations, is formalized in \S\ref{s.dist}.  

 In order to apply the John-Nirenberg Inequality, with the weight $ w$ fixed, we should decompose the collection of dyadic cubes into a \emph{Corona Decomposition.} As we work with Haar Shifts, a decomposition of the cubes leads automatically to the decomposition of the operator $ \mathbb S $.  This leads to a decomposition of $ \mathbb S (w \mathbf 1_{E})$ into terms which are individually very nicely behaved. 

 Finally, the testing conditions can be verified, and using the exponential integrability from the Corona Decomposition, one can give a simple verification of these conditions. This part of the argument is new to this paper.  
 This argument will \emph{not} appeal to extrapolation, a common technique in this subject. Indeed, the 
weak type estimate we prove does not seem to lend itself to extrapolation.

In the ultimate section, we provide some variations and consequences of Theorem~\ref{t.main}.

%

\section{Haar Shift Operators} \label{s.haarShift}

In this section, we introduce fundamental dyadic approximations of Calder\'on--Zygmund operators, the Haar shifts, and make a detailed reduction of the Main Theorem~\ref{t.main} to a similar statement, Theorem~\ref{t.haarShiftWtd},   in this dyadic model.

\begin{definition}\label{d.grid} A \emph{dyadic grid} is a collection $ \mathcal D$ of cubes so that for each $Q$ we have 
that 
\begin{enumerate}
\item  The set of cubes $ \{Q'  \in \mathcal D\;:\;  \lvert  Q'\rvert= \lvert  Q\rvert  \}$ partition $ \mathbb R ^{d}$, ignoring overlapping 
boundaries of cubes. 
\item $ Q$ is a  union of  cubes in a collection $ \textup{Child} (Q)  \subset \mathcal D$, called the \emph{children of $ Q$}. 
There are $ 2 ^{d}$ children of $ Q$, each of  volume $ \lvert  Q'\rvert= 2 ^{-d} \lvert  Q\rvert  $. 
\end{enumerate}
We refer to any subset of a dyadic grid as simply a \emph{grid}.  
\end{definition}

The standard choice for $ \mathcal D$ consists of the cubes $ 2 ^{k} \prod _{s=1} ^{d} [n_s, n_s+1)$ for $ k, n_1 ,\dotsc, n_d\in \mathbb Z $.  
But, the main result of this section, Theorem~\ref{t.bcr}, depends upon a random family of dyadic grids.

In higher dimensions, we mention that the martingale differences are finite rank projections, but there is no canonical choice of the Haar functions in this case.  We make the following definition.

\begin{definition}
Let $Q$ be a dyadic cube, a \emph{generalized Haar function} associated to $Q$, $h_{Q}$, is a linear combination of the indicator functions $ \displaystyle \left\{ \mathbf{1}_{Q}\right\} \cup \left\{ \mathbf{1}_{Q^{\prime
}}:Q^{\prime }\text{ is a child of }Q\right\}$, 
\begin{equation*}
h_{Q}=\sum_{Q^{\prime }\in \textup{Child}( Q ) }c_{Q^{\prime }}\mathbf{1}%
_{Q^{\prime }}
\end{equation*}%

We say $h_{Q}$ is a \emph{Haar function} if in addition $\int h_{Q}=0$, that is, a Haar function is
orthogonal to constants on its support. 
\end{definition}

\begin{definition}\label{d.haarShift} For integers $(m,n) \in \mathbb Z _+ ^2 $, we say that  a linear operator $ \mathbb S $ is a \emph{(generalized) Haar shift operator of   complexity type $ (m,n)$} if 
\begin{equation}\label{e.mn}
  \mathbb S  f (x) = \sum_{Q \in \mathcal D}\mathbb{S}_Q f(x)
  = \sum_{Q \in \mathcal D}\quad\sideset {} { ^ {(m,n)}} \sum_{\substack{Q',R'\in \mathcal D\\ Q',R'\subset Q }} 
\frac { \langle f, h ^{Q'} _{R'} \rangle} {\lvert  Q\rvert } k^{R'} _{Q'}(x)
\end{equation}
where  here and throughout $\ell (Q)= \lvert  Q\rvert ^{1/d} $, and   
\begin{itemize}
  \item in the second sum, the superscript $ ^{(m,n)}$ on the sum means that in addition we require $ \ell (Q') = 2 ^{-m} \ell (Q)$ and $ \ell (R')= 2 ^{-n} \ell (Q)$, and
  \item the function $  h ^{Q'} _{R'}$ is a (generalized) Haar function on $ R'$, and $  k ^{R'} _{Q'}$ is one on $ Q'$, with the normalization that 
\begin{equation} \label{e.normal}
\lVert  h ^{Q'} _{R'}\rVert_{\infty }\leq 1,\qquad \lVert  k ^{R'} _{Q'}\rVert_{\infty } \leq 1 \,.
\end{equation}
\end{itemize}
Here, and throughout the paper, $\ell (Q)=\lvert Q\rvert ^{1/d}$ is the side
length of the cube $Q$.
We say that the \emph{complexity} of $ \mathbb S $ is $  \kappa := \max (m,n,1)$.  
\end{definition}

A generalized Haar shift thus has the form%
\begin{equation*}
\mathbb{S}f\left( x\right) =\sum_{Q\in \mathcal{ D}}\frac{1}{\left\vert
Q\right\vert }\int_{Q}s_{Q}\left( x,y\right) f\left( y\right) dy=\int_{%
\mathbb{R}^{n}}K_{\mathbb S}\left(x,y\right) f\left( y\right) dy,
\end{equation*}%
where $s_{Q}$, the kernel of the component $\mathbb{S}_Q$, is supported on $Q\times Q$ and $\left\Vert s_{Q}\right\Vert
_{\infty }\leq 1$. It is easy to check that
\begin{equation*} 
  |K_{\mathbb{S}}(x,y)|\lesssim\frac{1}{|x-y|^d}.
\end{equation*}

The Haar shifts are automatically bounded on $L^2$ with $\|\mathbb{S}f\|_{L^2}\leq\|f\|_{L^2}$. This follows from the imposed normalizations and simple orthogonality considerations. It is clear that all \emph{restricted shifts} $\mathbb{S}_{\mathcal{Q}}f=\sum_{Q\in\mathcal{Q}}\mathbb{S}_Q f$ are also Haar shifts, and hence uniformly bounded on $L^2$ for any $\mathcal{Q}\subset\mathcal{D}$. We will extensively exploit these restricted shifts in the argument.

The generalized Haar shifts are only of relevance to us in two particular special cases of complexity type $(0,0)$, where $\mathbb{S}_Q f=|Q|^{-1}\langle f, h_Q^Q\rangle k_Q^Q=|Q|^{-1}\langle f, h_Q\rangle k_Q$. These are the \emph{paraproduct}, where $h_Q=\mathbf 1_Q$ and $k_Q$ is a Haar function for all $Q$, and the \emph{dual paraproduct}, where $h_Q$ is a Haar function and $k_Q=\mathbf 1_Q$ for all $Q$.

It is well-known that the (normalized) $L^2$ boundedness of a (dual) paraproduct is equivalent to the Carleson condition
\begin{equation*}
  \sum_{Q\subset R}\|k_Q\|_{L^2}^2\leq|R|\qquad\Big(\sum_{Q\subset R}\|h_Q\|_{L^2}^2\leq|R|\Big)\qquad\forall R\in\mathcal{D}.
\end{equation*}
These conditions are also uniformly inherited by all restricted (dual) paraproducts $\mathbb{S}_{\mathcal{Q}}$.

Note that for both Haar shifts and the paraproduct, we have
\begin{equation*}
  \int\mathbb{S}_Qf=0\qquad\forall Q\in\mathcal{D},
\end{equation*}
an important cancellation property in some of the later arguments. This is not the case for the dual paraproduct, for which a separate case study is needed at some points.


\begin{remark}
Let%
\begin{equation*}
\delta\left( x,y\right) =\min \left\{ \ell \left( Q\right) :x,y\in Q\in \mathcal{D%
}\right\}
\end{equation*}%
be the \emph{dyadic distance} between $x$ and $y$. The kernel $K_{\mathbb{S}}%
\left( x,y\right) $ of $\mathbb{S}$ satisfies the size and smoothness
conditions for a dyadic Calder\'on--Zygmund kernel:%
\begin{eqnarray*}
  \left\vert K_{\mathbb S}\left(x,y\right) \right\vert &\leq &\frac{2}{\delta\left(x,y\right) ^{d}}, \\
  \left\vert K_{\mathbb S}\left(x,y\right) -K_{\mathbb S}\left(x^{\prime},y\right) \right\vert &=&0\text{ if }
     \frac{\delta\left( x,x^{\prime }\right) }{\delta\left( x,y\right) }<\frac{1}{2^{m}}, \\
  \left\vert K_{\mathbb S}\left(x,y\right) -K_{\mathbb S}\left(x,y^{\prime}\right) \right\vert &=&0\text{ if }
      \frac{\delta\left( y,y^{\prime }\right) }{\delta\left( x,y\right) }<\frac{1}{2^{n}}.
\end{eqnarray*}%
This is a more general dyadic kernel condition than one studied in \cite{MR1934198}, called \emph{perfect dyadic}, which corresponds to $m=n=0$ in our framework.
\end{remark}

The relevance of Haar shifts to Classical Analysis is explained by the following Theorem, one of the main results of \cite{1010.0755} (see \cite{1010.0755}*{Theorem 4.1}; also \cite{1007.4330}*{Theorem 4.2}). 
  This Theorem must be formulated in terms of a \emph{random} dyadic grids.  But the nature of this construction of 
grids is immaterial to the arguments of this paper, and refer the reader to these references for proofs, history, and further discussion 
of this result.

\begin{theorem}\label{t.bcr}  
There is a collection of random dyadic grids $ \{\mathcal D _ \beta \;:\; \beta \in \boldsymbol \beta \}$, 
with expectation operator $ \mathbb E _{\beta }$, for which the following holds.  
Let $ T$ be a Calder\'on-Zygmund Operator $ T$ with smoothness parameter $ \delta $.  Then, for all bounded and compactly supported functions $f$ and $g$,  we can write 
\begin{equation} \label{e.bcr}
  \langle Tf,g\rangle = C  \mathbb E _{\beta } \sum_{ (m,n)\in \mathbb Z _+ ^2 }  2 ^{- (m+n) \delta /2} \langle \mathbb S  _{m,n} ^{\beta } f,g\rangle
\end{equation}
where
\begin{itemize}
  \item $\mathbb{S} _{m,n}^{\beta}$ is a Haar shift of complexity type $(m,n)$ for all $(m,n)\in\mathbb{Z}_+^2\setminus\{(0,0)\}$;
  \item $\mathbb{S}_{0,0}^{\beta}$ is the sum of a Haar shift of type $(0,0)$, a paraproduct, and a dual paraproduct;
  \item the constant $ C$ is a function of $ T$, and of the smoothness parameter $ \delta $. 
\end{itemize}
In particular we have the uniform estimate $\|\mathbb{S}_{m,n}^{\beta}\|_{L^2\to L^2}\leq 1$.
\end{theorem}

We define the maximal truncations of a Haar Shift as follows. 

\begin{definition}\label{d.maximalTruncations}
Suppose that $\mathbb{S}$ is a generalized Haar shift. Define the associated 
\emph{maximal truncations} by 
\begin{equation*}
\mathbb{S}_{\natural }f(x)\equiv \sup_{0<\epsilon\leq\upsilon <\infty }\lvert 
\mathbb{S}_{\epsilon ,\upsilon }f(x)\rvert \,,
\end{equation*}%
\begin{equation}
  \mathbb{S}_{\epsilon ,\upsilon }f\equiv \sum_{Q\in \mathcal{ D}:\ \epsilon
   \leq \ell (Q)\leq\upsilon }\mathbb{S}_Qf(x).  \label{e.Ter}
\end{equation}%
\end{definition}

\begin{proposition}\label{prop:TvsSpointwise}
We have the pointwise bound
\begin{equation*}
  T_{\natural}f(x) \le C(T)\mathbb{E} _{\beta}\sum_{(m,n) \in  \mathbb Z^2_+} 2^{-(m+n)\delta/2} (\mathbb S^{\beta}_{m,n})_{\natural}f(x)
   + C(T)Mf(x),
\end{equation*}
where $M$ is the Hardy-Littlewood maximal operator.
\end{proposition}

\begin{proof}
Theorem \ref{t.bcr} says that
\begin{displaymath}
\langle Tf, g\rangle = C(T)\mathbb{E} _{\beta} \sum_{(m,n) \in  \mathbb Z^2_+} 2^{-(m+n)\delta/2} \langle \mathbb S^{\beta}_{m,n}f, g\rangle
\end{displaymath}
for bounded and compactly supported functions $f$ and $g$. By choosing
\begin{displaymath}
  f = \frac{\mathbf{1}_{B(y, \epsilon)}}{|B(y,\epsilon)|} \qquad \textrm{and} \qquad g = \frac{\mathbf{1}_{B(x, \epsilon)}}{|B(x,\epsilon)|}
\end{displaymath}
and taking the limit as $\epsilon\to 0$, dominated convergence implies the pointwise identity
\begin{equation}\label{kernelidentity}
K(x,y) = C(T)\mathbb{E} _{\beta} \sum_{(m,n) \in  \mathbb Z^2_+} 2^{-(m+n)\delta/2} \sum_{Q \in \mathcal{ D}^{\beta}} \frac{s^{m,n}_Q(x,y)}{|Q|}.
\end{equation} 

For $\epsilon > 0$, this implies by Fubini's theorem that
\begin{align*}
\int_{|x-y| \ge \epsilon} &K(x,y)f(y)\,dy \\
& = C(T) \mathbb{E}_{\beta} \sum_{(m,n) \in  \mathbb Z^2_+} 2^{-(m+n)\delta/2} \sum_{Q \in \mathcal{ D}^{\beta}} \frac{1}{|Q|} \int_{B(x,\epsilon)^c} s_Q^{m,n}(x,y)f(y)\,dy.
\end{align*}
Let us then decompose
\begin{align*}
\sum_{Q \in \mathcal{ D}^{\beta}} \frac{1}{|Q|} \int_{B(x,\epsilon)^c} s_Q^{m,n}(x,y)f(y)\,dy &= \mathop{\sum_{Q \in \mathcal{ D}^{\beta}}}_{\ell(Q) > \epsilon'} \frac{1}{|Q|} \int s_Q^{m,n}(x,y)f(y)\,dy \\
& \quad -  \mathop{\sum_{Q \in \mathcal{ D}^{\beta}}}_{\ell(Q) > \epsilon'} \frac{1}{|Q|} \int_{B(x,\epsilon)} s_Q^{m,n}(x,y)f(y)\,dy \\
&\quad + \mathop{\sum_{Q \in \mathcal{ D}^{\beta}}}_{\ell(Q) \le \epsilon'} \frac{1}{|Q|} \int_{B(x,\epsilon)^c} s_Q^{m,n}(x,y)f(y)\,dy \\
&= I  +I\!I + I\!I\!I,
\end{align*}
where $\epsilon' = \epsilon/(2\sqrt{d})$. By definition, there holds $|I| \le (\mathbb{S}^{\beta}_{m,n})_{\natural}f(x)$. There also holds
\begin{align*}
|I\!I| \le \sum_{k:\, 2^k > \epsilon'} 2^{-kd}\int_{B(x,\epsilon)} |f| \lesssim \epsilon^d \Big( \sum_{k:\, 2^k > \epsilon'} 2^{-kd} \Big)Mf(x) \lesssim Mf(x).
\end{align*}
Furthermore, we actually have $I\!I\!I = 0$, since there we must have $|x-y| \le d(Q) = \sqrt{d}\ell(Q) \le \sqrt{d}\epsilon' = \epsilon/2 < \epsilon$.
We have thus shown that
\begin{displaymath}
T_{\epsilon}f(x) \le C(T)\mathbb{E} _{\beta} \sum_{(m,n) \in  \mathbb Z^2_+} 2^{-(m+n)\delta/2} (\mathbb S^{\beta}_{m,n})_{\natural}f(x) + C(T)Mf(x)
\end{displaymath}
for every $\epsilon > 0$, and from this the proposition follows.
\end{proof}

Proposition~\ref{prop:TvsSpointwise} and Buckley's \cite{MR1124164} sharp weighted bounds for the maximal operator,
\begin{equation}\label{eq:Buckley}
  \|Mf\|_{L^{p,\infty}(w)}\lesssim\|w\|_{A_p}^{1/p}\|f\|_{L^p(w)},\qquad
  \|Mf\|_{L^{p}(w)}\lesssim\|w\|_{A_p}^{1/(p-1)}\|f\|_{L^p(w)},
\end{equation}
reduce the proof of the Main Theorem~\ref{t.main} to the verification of the following dyadic variant, a task which occupies the rest of this paper.


\begin{theorem}\label{t.haarShiftWtd}  Let $ \mathbb S $ be a Haar shift operator with complexity $ \kappa $, a paraproduct, or a dual paraproduct.
For $ 1< p < \infty $ and $ w \in A_p$,  we then have the estimates 
\begin{align}\label{e.Sweak}
\lVert \mathbb S_{\natural } f \rVert_{ L ^{p, \infty } (w)} &\lesssim \kappa \lVert w\rVert_{A_p}\lVert f\rVert_{L ^{p} (w)}  \,, 
\\ \label{e.Sstrong}
\lVert \mathbb S_{\natural } f \rVert_{ L ^{p } (w)} & \lesssim \kappa \lVert w\rVert_{A_p}  ^{\max \{1, (p-1) ^{-1} \}}\lVert f\rVert_{L ^{p} (w)} 
\end{align}
\end{theorem}

Indeed, any polynomial dependence on the complexity parameter $\kappa$ would suffice for Theorem~\ref{t.main}, but a careful tracing of the constants will even provide the \emph{linear} dependence, as stated. Even for the untruncated shifts $\mathbb{S}$ in $L^2(w)$, this improves on the quadratic in $\kappa$ bound established in \cite{1010.0755} (but we are not aware of an application where this precision in the dependence on $\kappa$ would be of importance).

The dependence on $ \kappa$ and the weight constant $\|w\|_{A_p}$ arises from the following points of the proof below: First, we establish \emph{two-weight inequalities} of the form
\begin{align}
  \lVert \mathbb S_{\natural }(f\sigma) \rVert_{ L ^{p, \infty } (w)}
  &\lesssim \left\{\kappa\mathfrak{M}_{p,\operatorname{weak}}+\mathfrak{T}_p\right\}\lVert f\rVert_{L ^{p} (\sigma)}, \\
  \lVert \mathbb S_{\natural } (f\sigma) \rVert_{ L ^{p } (w)}
  &\lesssim \left\{\kappa\mathfrak{M}_{p}+\mathfrak{T}_p+\mathfrak{N}_p\right\}\lVert f\rVert_{L ^{p} (\sigma)},
\end{align}
where $\mathfrak{M}_{p,\operatorname{weak}}$ and $\mathfrak{M}_p$ are the best constants from certain maximal inequalities, while $\mathfrak{T}_p$ and $\mathfrak{N}_p$ are the best constants from appropriate \emph{testing conditions} for the operator $\mathbb{S}_{\natural}$. Here, $w$ and $\sigma$ are allowed to be an arbitrary pair of weights, with no relation to each other.
 
Second, we specialize to the one-weight situation with $\sigma=w^{1-p'}$, using a well-known dual-weight formulation of the bounds to be proven, \eqref{e.Sweak} and \eqref{e.Sstrong}.
We need to estimate the above four constants in this situation. The maximal constants are independent of $\mathbb{S}_{\natural}$ and thus of $\kappa$, and they satisfy $\mathfrak{M}_{p,\operatorname{weak}}\lesssim\|w\|_{A_p}^{1/p}$ and $\mathfrak{M}_p\lesssim\|w\|_{A_p}^{1/(p-1)}$ by the sharp maximal function inequalities of Buckley \eqref{eq:Buckley}. For the two testing constants related to $\mathbb{S}_{\natural}$, we obtain the linear in $\kappa$ bounds
\begin{equation*}
  \mathfrak{T}_p\lesssim\kappa\|w\|_{A_p},\qquad\mathfrak{N}_p\lesssim\kappa\|w\|_{A_p}^{1/(p-1)}.
\end{equation*}
This dependence comes from the fact that the proof of the John--Nirenberg style estimates of \eqref{e.Lexp2} requires \emph{separating the scales} of $\mathbb{S}$ by dividing it into $\kappa+1$ parts, each of which contains nonzero components $\mathbb{S}_Q$ only for a fixed value of $\log_2\ell(Q)\mod(\kappa+1)$. For these separated parts of $\mathbb{S}$, our bounds will be independent of $\kappa$, and it remains to sum up.

\section{Linearizing Maximal Operators}

\label{s.linear}

A fundamental tool is derived from (the usual) general maximal function
estimates that hold for any measure. In particular, for weight $w $ we
define 
\begin{align}
M_{w }f(x)& \equiv \sup_{Q\in \mathcal{D}}\mathbf{1}_{Q}(x)\mathbb{E}%
_{Q}^{w }\lvert f\rvert \,,  \label{e.wtdE} \\
\mathbb{E}_{Q}^{w }f& \equiv w (Q)^{-1}\int_{Q}f\,w (dx)\,.
\end{align}%
 The notation $ \mathbb E _{Q} f$ means that the implied measure is Lebesgue.
It is a basic fact, proved by exactly the same methods that prove
the non-weighted inequality, that we have the estimate below, which will be used repeatedly. 


\begin{theorem}
\label{t.max} We have the inequalities 
\begin{equation}
\left\Vert M_{w }f\right\Vert _{L^{p}(w )}\lesssim \left\Vert
f\right\Vert _{L^{p}(w )}\,,\qquad 1<p<\infty \,.  \label{e.Umax}
\end{equation}
\end{theorem}


We use the method of linearizing maximal operators. This is
familiar in the context of the maximal function, and we make a comment about
it here. Let $\{E(Q)\;:\;Q\in \mathcal{D}\}$ be any selection of measurable
disjoint sets $E(Q)\subset Q$ indexed by the dyadic cubes. Define a
corresponding linear operator $N$ by 
\begin{equation}
N\phi \equiv \sum_{Q\in \mathcal{D}}\mathbf{1}_{E(Q)}\mathbb{E}_{Q}^{w
}f\,.  \label{e.N}
\end{equation}%
Then, the universal Maximal function bound \eqref{e.Umax} is equivalent to
the bound $\left\Vert Nf\right\Vert _{L^{p}(w )}\lesssim \left\Vert
f\right\Vert _{L^{p}(w )}$ with implied constant independent of $w 
$ and the sets $\{E(Q)\;:\;Q\in \mathcal{D}\}$. This estimate will be used
repeatedly below.

\medskip

There is a related way to linearize $\mathbb{S}_{\natural }$, which deserves
careful comment as we would like, at different points, to treat $\mathbb{S}%
_{\natural }$ as a linear operator. While it is not a linear operator, $%
\mathbb{S}_{\natural }$ is a pointwise supremum of the linear truncation
operators $\mathbb{S}_{\varepsilon ,\upsilon }$, and as such, the supremum
can be linearized with measurable selection of the truncation parameters.

\begin{definition}
\label{d.linearization} We say that $\mathbb{L}$ is a linearization of $%
\mathbb{S}_{\natural }$ if there are measurable functions $\epsilon
(x),\upsilon (x)\in (0,\infty )$ and $\vartheta (x)\in \lbrack 0,2\pi )$
such that, using \eqref{e.Ter}, we have 
\begin{equation}
\mathbb{L}f(x)=e^{i\vartheta (x)}\mathbb{S}_{\epsilon (x),\upsilon
(x)}f(x)\geq 0,\qquad x\in \mathbb{R}^{d}.  \label{defLa}
\end{equation}
\end{definition}

Note that the requirement $\mathbb{L}f(x)\geq 0$ defines $\vartheta (x)$
everywhere except when $\mathbb{S}_{\epsilon (x),\upsilon (x)}f(x)=0$. Also,
for fixed $f$ we can always choose a linearization $L$ so that $\mathbb{S}%
_{\natural }f(x)\leq 2\mathbb{L}f(x)$ for all $x$.

A key advantage of $\mathbb{L}$ is that it is a linear operator, and as such
it has an adjoint, given by the formal expressions
\begin{equation}\label{e.adjoint}
\begin{split}
  \mathbb{L}^*\nu(y)
  &=\sum_{Q\in\mathcal{D}}\mathbb{S}_Q^*(\mathbf 1_{\{\epsilon(\cdot)\leq\ell(Q)\leq\upsilon(\cdot)\}}e^{i\vartheta(\cdot)}\nu)(y) \\
  &=\sum_{Q\in\mathcal{D}}\frac{1}{|Q|}\int_Q s_Q(x,y)\mathbf 1_{\{\epsilon(x)\leq\ell(Q)\leq\upsilon(x)\}}e^{i\vartheta(x)}\nu(dx).
\end{split}
\end{equation}
%

The following `smoothness' property of $ \mathbb L ^{\ast} $ is an important observation in the proof 
of our two weight estimates.

\begin{lemma}\label{lem:smoothness}
\label{l.L} Suppose that for a measure $\nu $ and cube $Q_{0}$ we have $\lvert
\nu \rvert (Q_{0})=0$. Suppose that $\mathbb{S}$ has complexity type $\left(
m,n\right) $. Then $\mathbb{L}^{\ast }\nu (\cdot )$ is constant on subcubes $%
Q^{\prime }\subset Q_{0}$ with $\ell(Q')\leq 2^{-n}\ell(Q_0)$.
\end{lemma}

\begin{proof}
For $y\in Q_0$, the sum in \eqref{e.adjoint} defining the adjoint operator becomes 
\begin{equation*}
  \mathbb{L}^*\nu(y)
  =\sum_{Q\in\mathcal{D}:\ Q\supsetneq Q_0}\frac{1}{|Q|}\int_Q s_Q(x,y)\mathbf 1_{\{\epsilon(x)\leq\ell(Q)\leq\upsilon(x)\}}e^{i\vartheta(x)}\nu(dx).
\end{equation*}
As a function of $y$, the kernel $s_Q(x,y)$ is constant on the subcubes of $R'\subset Q$ with $\ell(R')\leq 2^{-n-1}\ell(Q)$. Thus for $Q\supsetneq Q_0$, it is in particular constant on the subcubes $Q'\subset Q_0$ with $\ell(Q')\leq 2^{-n}\ell(Q_0)$.
\end{proof}

\section{The Two Weight Estimates} \label{s.2wt}

We are interested in tracking complexity dependence in two weight
inequalities for Haar Shift Operators, as defined in \S\ref{s.haarShift}. We study the maximal
truncations of such operators, and obtain sufficient conditions for the weak
and strong type $(p,p)$ two weight inequalities for such operators. Our main
results are Theorem~\ref{weakFlat} for the weak type result, and Theorem \ref%
{t.genlStrong} for the strong type result.   
These Theorems give sufficient conditions in terms of the Maximal Function, and certain 
\emph{testing conditions}.   Of particular import here is that these sufficient conditions are efficient in terms of 
the \emph{complexity} of the Haar Shift operator.  


Our primary focus concerns extensions of the dyadic $T1$ Theorem to the two
weight setting. These considerations are motivated in part by a well
developed theory of two weight estimates for positive operators. These
Theorems have formulations strikingly similar to the $T1$ Theorem, which
theory encompasses the Theorems due to one of us concerning two weight, both
strong and weak type, for the maximal operator \cite{MR676801} and
fractional integral operators \cite{MR719674}, \cite{MR930072}. There is
also the bilinear embedding inequality of Nazarov-Treil-Volberg \cite{MR1685781}.
We refer the reader to \cite{0911.3437} for a discussion of these results.

There is a beautiful result of Nazarov-Treil-Volberg \cite{NTV2}, a
two-weight version of the $T1$ theorem.
 A subcase of their result was proved for Haar Shifts, with an effective 
bound on complexity in \cite{1010.0755}*{Theorem 3.4}. 

\begin{theorem}
\label{t.ntv} Let $\mathbb{S}$ be a Haar Shift operator of complexity $ \kappa $, as in Definition~\ref{d.haarShift}. Let $\sigma
,w $ be two positive locally finite measures. We have the inequality 
\begin{equation}
\lVert \mathbb S  (f \sigma )\rVert_{ L ^2 (w )} 
\lesssim \kappa \{ \mathfrak  S + \mathfrak  S ^{\ast} \} + \kappa ^2 \lVert w , \sigma \rVert_{A_2}^{{1/2}}
\label{e.ntv1}
\end{equation}%
where the three quantities above are defined by 
\begin{align}
\lVert w, \sigma \rVert_{A_2} 
& := \sup _{Q} \frac { w (Q)} {\lvert  Q\rvert } \frac {\sigma (Q)} {\lvert  Q\rvert}
\\  
\mathfrak S &:= 
\sup _{Q} \sigma  (Q) ^{-1/2} \lVert \mathbf 1_{Q} \mathbb S (\sigma\mathbf 1_{Q})\rVert_{ L ^2 (w)} \,. 
\\ \label{e.ntv4}
\mathfrak S ^{\ast}  &:= 
\sup _{Q} w (Q) ^{-1/2} \lVert \mathbf 1_{Q} \mathbb S ^{\ast}  (w \mathbf 1_{Q})\rVert_{ L ^2 (\sigma )} \,. 
\end{align}
\end{theorem}

The first of the three conditions is the two weight $ A_2$ condition; the remaining two are the testing conditions. 
The proof is fundamentally restricted to the case of $p=2$, nor does it address maximal truncations. 
We will consider the case of $1<p<\infty $ and obtain  sufficient conditions for    the two weight
inequalities for the maximal operator $\mathbb{S}_{\natural }$. First we
give the weak type result. Below, $ M$ denotes the maximal function.  

\begin{theorem}
\label{weakFlat} Let $\mathbb{S}$ be a generalized Haar shift of complexity $%
\kappa $ as in Definition~\ref{d.haarShift}. Then we have the weak type inequality 
\begin{equation}
\left\Vert \mathbb{S}_{\natural }(f\sigma )\right\Vert _{L^{p,\infty}(w )}
\lesssim \left(\kappa \mathfrak{M}_{p,\textup{weak}}+\mathfrak{T}_p\right) \left\Vert f\right\Vert _{L^{p}(\sigma )}\,,
\label{e.TW2}
\end{equation}%
where the constants $\mathfrak{M}_{p,\textup{weak}}$ and $\mathfrak{T}_p$
are the best such in the following inequalities%
\begin{align}
\left\Vert {M}(f\sigma )\right\Vert _{L^{p,\infty }(w )}\leq \mathfrak{M%
}_{p,\textup{weak}}\left\Vert f\right\Vert _{L^{p}(\sigma )}\,,  \label{e.Mweak}
\\
\int_{Q}\mathbb{S}_{\natural }(\sigma f\mathbf{1}_{Q})\,\;w (dx)\leq 
\mathfrak{T}_p\left\Vert f\right\Vert _{L^{p}(\sigma )}w
(Q)^{1/p^{\prime }}.  \label{e.tw2}
\end{align}
\end{theorem}

The point of this Theorem is that to check the weak-type inequality for $%
\mathbb{S}_{\natural }$, it suffices to check the weak-type inequality for
the simpler maximal operator $M$, and to check only particular instances of
the weak-type inequality for $\mathbb{S}_{\natural }$. It is also important
that the complexity $\kappa $ appears with polynomial growth.

The dual testing condition \eqref{e.tw2} looks rather complicated, with the appearance of $f\in
L^{p}(\sigma )$ in it. However, $\mathbb{S}_{\natural }$ appears to just the
first power, and it is a close relative of \eqref{e.ntv4}.
Indeed, \eqref{e.tw2} has a more convincing formulation in
the linearizations. It is equivalent to the dual testing condition%
\begin{equation}
  \left\Vert \mathbf{1}_Q\mathbb{L}^{\ast }(\mathbf{1}_{Q}gw )\right\Vert_{L^{p^{\prime }}(\sigma )}
  \leq \mathfrak{T}_p w (Q)^{1/p^{\prime}}\left\Vert g\right\Vert _{\infty },  \label{e.wf3'}
\end{equation}%
This holds uniformly over all choices of linearizations, which fact is
referred to repeatedly below. 
Inequality
\eqref{e.wf3'} reflects the fact that the dual of a weak type inequality is
a restricted strong type inequality.

Our strong type result will require duals $\mathbb{L}^{\ast }$ of
linearizations $\mathbb{L}$ of $\mathbb{S}_{\natural }$ in order to state
the nonstandard testing condition in \eqref{e.non}. These were defined in Section \ref{s.linear} above. 

\begin{theorem}
\label{t.genlStrong} Let $\mathbb{S}$ be a generalized Haar shift of
complexity $\kappa $ as in Definition~\ref{d.haarShift}. We have the following
quantitative estimate: 
\begin{equation}
\lVert \mathbb{S}_{\natural }(f\sigma )\rVert _{L^{p}(w )}\lesssim
\left\{ \kappa\mathfrak{M}_{p}+\mathfrak{T}_p+\mathfrak{N}_{p}\right\} \lVert f\rVert _{L^{p}(\sigma )}\,,  \label{e.STR}
\end{equation}%
where $\mathfrak{T}_p$ is defined in \eqref{e.tw2}, and the numbers $\mathfrak{M}_{p}$ and $\mathfrak{N}_{p}$ are defined in 
\begin{equation}
\left\Vert {M}(f\sigma )\right\Vert _{L^{p}(w )}\leq \mathfrak{M}%
_{p}\left\Vert f\right\Vert _{L^{p}(\sigma )},  \label{e.Mstrong}
\end{equation}%
\begin{equation}
\mathfrak{N}_{p}^{p}\equiv \sup_{\left\Vert \varphi \right\Vert _{\infty
}\leq 1}\sup_{Q_{0}}\frac{1}{\sigma (Q_{0})}\int_{Q_{0}}\sup_{Q\subset Q_{0}}%
\mathbf{1}_{Q}\left( \frac{1}{w (Q)}\int_{Q}\left\vert \mathbb{L}^{\ast
}(\mathbf{1}_{Q}\varphi w )(y)\right\vert \sigma (dy)\right) ^{p}w
(dx).  \label{e.non}
\end{equation}
\end{theorem}


As a new kind of complication compared to the weak-type case, we have the nonstandard testing condition \eqref{e.non}.
Its primary difficulty is the appearance of $\left\vert
\mathbb L^{\ast }\right\vert $ integrated over $Q$ with respect to $\sigma $, but
then divided by $w (Q)$ rather than the usual $\sigma (Q)$. Also there
is an additional supremum, with the argument of $\mathbb L^{\ast }$ \emph{dependent}
upon the cube $Q$ over which we are taking the supremum.

The method of proof is an extension of that of Sawyer's approach to the two
weight fractional integrals \cite{MR930072}, but see also \cite{0911.3437}.
This argument follows
the outlines of the proof in \cite{0807.0246}, which proves variants of
Theorem~\ref{weakFlat} and Theorem~\ref{t.genlStrong} for smooth Calder\'{o}%
n-Zygmund operators. The current arguments are, naturally, much easier while
retaining the essential ideas and techniques of \cite{0807.0246}. (The
reader can also compare the arguments of this paper to those of \cite%
{0911.3437}.)
 
Let us give a guide to the {next few sections of this paper, which are concerned with the proof of the above two-weight results.}

\begin{description}
\item[\S \protect\ref{s.general}] Collects facts central to the proofs,
maximal functions, linearizations of maximal functions, Whitney
decompositions, and an important maximum principle.

\item[\S \protect\ref{s.weak}] The weak-type result Theorem~\ref{weakFlat}
is proved.

\item[\S \protect\ref{s.first}] Sufficient conditions for the strong type
result are stated; the classical part of the proof of the strong type result
Theorem~\ref{t.genlStrong} is begun.

\item[\S \protect\ref{s.nondoubling}] This section contains
the core of the proof of Theorem\ \ref{t.genlStrong}.
\end{description}

\section{Generalities of the Proof}\label{s.general}

\subsection{Whitney Decompositions}

We make general remarks about the sets 
\begin{equation}
\Omega _{k}=\{\mathbb{S}_{\natural }(f\sigma )>2^{k}\}  \label{e.Wdef}
\end{equation}%
where $f$ is a finite linear combination of indicators of dyadic cubes. 
For points $ x$ sufficiently far away from the support of $ f$, we will have that $ \mathbb S _{\natural} (f \sigma )$ 
is dominated by the maximal function $ M (f \sigma )$.  Hence, the sets $ \Omega _k $
will be  open  with compact closure.

Let $Q^{(1)}$ denote the dyadic parent of $Q$, and inductively define $Q^{(j+1)}=(Q^{(j)})^{(1)}$.
For a nonnegative integer $\zeta$, let $\mathcal{Q}_k$ be the collection of maximal dyadic cubes $Q$ such that $Q^{(\zeta)}\subset\Omega_k$. Then
\begin{equation}
\Omega _{k}=\dot{\bigcup_{Q\in \mathcal{Q}_{k}}}Q\ \ \ \ \ \text{\textup{%
(disjoint cover),}}  \label{e.dc}
\end{equation}%
\begin{equation}
Q^{(\zeta )}\subset \Omega _{k}\,,\ Q^{(\zeta +1)}\cap \Omega _{k}^{c}\neq
\emptyset \ \ \ \ \ \text{\textup{(Whitney condition,)}}  \label{e.Whit}
\end{equation}%
\begin{equation}
Q\in \mathcal{Q}_{k}\,,\ Q^{\prime }\in \mathcal{Q}_{l}\,,\ Q\subsetneqq
Q^{\prime }\quad \textup{implies}\quad k>l\ \ \ \ \ \text{\textup{(nested
property).}}  \label{e.nested}
\end{equation}

\begin{remark}
In the proof of the weak type theorem we will take $\zeta=0 $.
In the proof of the strong type theorem we
will take $\zeta=n+1\leq\kappa+1$, when the shift under consideration has complexity type $(m,n)$.
\end{remark}

\subsection{Maximum Principle}

A fundamental tool is the use of what we term here as a `maximum principle'
(we could also use the term `good-$\lambda $ technique'): Subject to the
assumption that the maximal function $M$ is of small size, we will be able
to see that the maximal truncations are large due to the restriction of the
function to a local cube. This leads to an essential `localization' of the
singular integrals.

\begin{theorem}[Maximum Principle]
Let $\mathbb{S}$ be a generalized Haar shift of complexity type $\left(
m,n\right) $. For any cube $Q\in \mathcal{Q}_{k}$ as in the Whitney
decomposition of $\Omega _{k}$ in \eqref{e.Wdef} above with parameter $\zeta$, we have the pointwise inequality 
\begin{equation}\label{maxprinc}
\begin{split}
  \mathbb{S}_{\natural }\left( f\sigma \right) \left( x\right)
  &\leq\sup_{\epsilon\leq\upsilon\leq\ell(Q)}|\mathbb{S}_{\epsilon,\upsilon}(f\sigma)(x)|+2^k+(\zeta+n+1)M(f\sigma)(x) \\
  &\leq\mathbb{S}_{\natural}(f1_Q\sigma)(x)+2^k+(\zeta+n+1)M(f\sigma)(x).
\end{split}
\end{equation}
\end{theorem}

There is a corresponding Maximum Principle in {Section 3.3 of} \cite%
{0911.3437}, which is very effective in the positive operator case. As our
operators are not positive, and as we are ultimately only interested in the
one weight situation, we have the maximal function $M$ on the right in %
\eqref{maxprinc}.

\begin{proof}
Note that the second inequality of the claim is obvious, since $\mathbb{S}_{\epsilon,\upsilon}(f\sigma)(x)=\mathbb{S}_{\epsilon,\upsilon}(f1_Q\sigma)(x)$ when $x\in Q$ and $\epsilon\leq\upsilon\leq\ell(Q)$. We prove the first inequality.

Let $x\in Q\in\mathcal{Q}_k$, and let $r:=(\zeta+1)+(n+1)$. Then
\begin{equation*}
\begin{split}
  \mathbb{S}_{\epsilon,\upsilon}(f\sigma)(x)
  &=\sum_{I:\epsilon\leq\ell(I)\leq\upsilon}\mathbb{S}_I(f\sigma)(x) \\
  &=\Big(\sum_{\substack{I:\epsilon\leq\ell(I)\leq\upsilon \\ \ell(I)\leq \ell(Q) }}
    +\sum_{\substack{I:\epsilon\leq\ell(I)\leq\upsilon \\ \ell(I)\geq 2^{r}\ell(Q) }}
    +\sum_{\substack{I:\epsilon\leq\ell(I)\leq\upsilon \\ \ell(Q)< \ell(I)< 2^{r}\ell(Q) }}\Big)\mathbb{S}_I(f\sigma)(x) \\
 &=\mathbb{S}_{\epsilon,\min(\upsilon,\ell(Q))}(f\sigma)(x)
    +\mathbb{S}_{\max(\epsilon,2^r\ell(Q)),\upsilon}(f\sigma)(x)
    +\sum_{\substack{I:\epsilon\leq\ell(I)\leq\upsilon \\ \ell(Q)< \ell(I)< 2^{r}\ell(Q) }}\mathbb{S}_I(f\sigma)(x).
\end{split}
\end{equation*}
The first term on the right is clearly dominated by $\sup_{\epsilon'\leq\upsilon'\leq\ell(Q)}|\mathbb{S}_{\epsilon',\upsilon'}(f\sigma)(x)|$. All $\mathbb{S}_I(f\sigma)$ participating in the second term have $\ell(I)\geq 2^{n+1}\ell(Q^{(\zeta+1)})$, so they are constant on $Q^{(\zeta+1)}$. Hence we may replace $x$ by some $\bar{x}\in Q^{(\zeta+1)}\setminus\Omega_k$ (which is nonempty by definition of $\mathcal{Q}_k$). So this second term is dominated by
\begin{equation*}
  \mathbb{S}_{\natural}(f\sigma)(\bar{x})\leq 2^k.
\end{equation*}
Finally, the last term contains at most $r-1=\zeta+n+1$ summands, each of which is dominated by $M(f\sigma)(x)$.
\end{proof}

\section{Proof of the Weak-Type Inequality}\label{s.weak}

We prove Theorem \ref{weakFlat}, stating that
\begin{equation*}
  \|\mathbb{S}_{\natural}(f\sigma)\|_{L^{p,\infty}(w)}
  \lesssim(\kappa\mathfrak{M}_{p,\textup{weak}}+\mathfrak{T}_p)\|f\|_{L^p(\sigma)}.
\end{equation*}
To this end, we need to estimate the quantities
\begin{equation*}
  w(\mathbb{S}_{\natural}(f\sigma)>4\cdot 2^k)
  \leq w(M(f\sigma)>\eta 2^k)+w(\mathbb{S}_{\natural}(f\sigma)>4\cdot 2^k,M(f\sigma)\leq\eta 2^k),
\end{equation*}
where the small parameter $\eta$ is to be chosen shortly. Since $\Omega_{k+2}\subset\Omega_k=\bigcup_{Q\in\mathcal{Q}_k}Q$, where we take the Whitney decomposition with parameter $\zeta=0$, we further have that
\begin{equation*}
\begin{split}
  w(\mathbb{S}_{\natural}(f\sigma) &>4\cdot 2^k,M(f\sigma)\leq\eta 2^k) \\
  &=\sum_{Q\in\mathcal{Q}_k}w(Q\cap\{\mathbb{S}_{\natural}(f\sigma)>4\cdot 2^k,M(f\sigma)\leq\eta 2^k\})
  =:\sum_{Q\in\mathcal{Q}_k}w(E_k(Q)).
\end{split}
\end{equation*}
On $E_k(Q)\subset Q$, the maximum principle gives that
\begin{equation*}
\begin{split}
  4\cdot 2^k<\mathbb{S}_{\natural}(f\sigma)
  &\leq\mathbb{S}_{\natural}(f 1_Q\sigma)+2^k+(2n+3)M(f\sigma) \\
  &\leq\mathbb{S}_{\natural}(f 1_Q\sigma)+2^k+(2n+3)\eta 2^k
    \leq \mathbb{S}_{\natural}(f 1_Q\sigma)+2\cdot 2^k,
\end{split}
\end{equation*}
provided that $\eta:=(2\kappa+3)^{-1}\leq (2n+3)^{-1}$. Thus
\begin{equation*}
  \mathbb{S}_{\natural}(f\sigma 1_Q)>2^{k+1}\qquad\text{on }E_k(Q).
\end{equation*}

Putting these considerations together, we obtain (for another small parameter $\delta>0$)
\begin{equation*}
\begin{split}
   (4\cdot 2^k)^p & w(\mathbb{S}_{\natural}(f\sigma)>4\cdot 2^k) \\
   &\leq 4^p 2^{kp}w(M(f\sigma)>\eta 2^k)+2^p\sum_{Q\in\mathcal{Q}_k} 2^{(k+1)p} w(E_k(Q)) \\
   &\leq 4^p\eta^{-p} \|M(f\sigma)\|_{L^{p,\infty}(w)}^p+4^p 2^{kp}\sum_{\substack{Q\in\mathcal{Q}_k \\ w(E_k(Q))<\delta w(Q)}} \delta w(Q) \\
   &\qquad     +2^p \sum_{\substack{Q\in\mathcal{Q}_k \\ w(E_k(Q))\geq \delta w(Q)}}
      w(E_k(Q))\Big(\frac{1}{w(E_k(Q))}\int_{E_k(Q)}\mathbb{S}_{\natural}(f 1_Q\sigma) w\Big)^p \\
    &\leq 2^p\eta^{-p}\mathfrak{M}_{p,\textup{weak}}^p\|f\|_{L^p(\sigma)}^p +4^p\delta 2^{kp} w(\mathbb{S}_{\natural}(f\sigma)>2^k) \\
    &\qquad +2^p \sum_{Q\in\mathcal{Q}_k}\delta^{1-p} w(Q)^{1-p}\Big(\int_Q\mathbb{S}_{\natural}(f1_Q\sigma)w\Big)^{p}
\end{split}
\end{equation*}
Picking a $k$ for which the left side is close to its supremum, and choosing $\delta=10^{-p}$,  we can absorb the middle term on the right to the left side.
Recalling the choice of $\eta$ and using the testing condition to the integrals in the last term, we obtain
\begin{equation*}
\begin{split}
  \|\mathbb{S}_{\natural}(f\sigma)\|_{L^{p,\infty}(w)}^p
  &\lesssim\kappa^p\mathfrak{M}_{p,\textup{weak}}^p\|f\|_{L^p(\sigma)}+\sum_{Q\in\mathcal{Q}_k} w(Q)^{1-p}(\mathfrak{T}_p\|f 1_Q\|_{L^p(\sigma)}w(Q)^{1/p'}\big)^p \\
  &=\kappa^p\mathfrak{M}_{p,\textup{weak}}^p\|f\|_{L^p(\sigma)}+\mathfrak{T}_p^p\sum_{Q\in\mathcal{Q}_k} \|f 1_Q\|_{L^p(\sigma)}^p \\
  &\leq\big(\kappa^p\mathfrak{M}_{p,\textup{weak}}^p+\mathfrak{T}_p^p\big)\|f\|_{L^p(\sigma)},
\end{split}
\end{equation*}
by the disjointness of the $Q\in\mathcal{Q}_k$ in the last step.

\section{First Steps in the Proof of the Strong Type Inequality}\label{s.first}

We start preparing for the proof of
\begin{equation*}
  \|\mathbb{S}_{\natural}(f\sigma)\|_{L^p(w)}
 \lesssim(\kappa\mathfrak{M}_p+\mathfrak{T}_p+\mathfrak{N}_p)\|f\|_{L^p(\sigma)}.
\end{equation*}
In this section, we make an estimate of the form
\begin{equation*}
  \|\mathbb{S}_{\natural}(f\sigma)\|_{L^p(w)}
  \lesssim \kappa\mathfrak{M}_p\|f\|_{L^p(\sigma)}+\delta\|\mathbb{S}_{\natural}(f\sigma)\|_{L^p(w)}
  +\text{remainder},
\end{equation*}
where the second term with the small parameter $\delta$ may be absorbed to the right, and the `remainder' will be controlled in terms on $(\mathfrak{T}_p+\mathfrak{N}_p)\|f\|_{L^p(\sigma)}$ in the following section, which contains the core of the argument.

We begin with
\begin{equation*}
\begin{split}
  \|\mathbb{S}_{\natural}(f\sigma)\|_{L^p(w)}^p
  &=\sum_{k\in\mathbb{Z}}\int_{\Omega_{k+1}\setminus\Omega_{k+2}}\mathbb{S}_{\natural}(f\sigma)^p w \\
  &\leq 4^p\sum_{k\in\mathbb{Z}} 2^{kp} w(\Omega_{k+1}\setminus\Omega_{k+2}) \\
  &=4^p\sum_{k\in\mathbb{Z}} 2^{kp} \sum_{Q\in\mathcal{Q}_k} w(Q\cap(\Omega_{k+1}\setminus\Omega_{k+2})) \\
  &=:4^p\sum_{k\in\mathbb{Z}} 2^{kp} \sum_{Q\in\mathcal{Q}_k} w(F_k(Q)).
\end{split}
\end{equation*}
Note that the sets $F_k(Q)$ are pairwise disjoint.
We first employ a similar reduction as in the weak-type case:
\begin{equation}\label{e.EQdef}
\begin{split}
  w(F_k(Q))
  &\leq w(F_k(Q)\cap\{M(f\sigma)>\eta 2^k\})+w(F_k(Q)\cap\{M(f\sigma)\leq\eta 2^k\}) \\
  &=: w(F_k(Q)\cap\{M(f\sigma)>\eta 2^k\})+w(E_k(Q)).
\end{split}
\end{equation}
Then
\begin{equation*}
\begin{split}
  \sum_{k} 2^{kp} \sum_{Q\in\mathcal{Q}_k} w(F_k(Q)\cap\{M(f\sigma)>\eta 2^k\})
  &\leq \eta^{-p}\sum_{k}\sum_{Q\in\mathcal{Q}_k}\int_{F_k(Q)} M(f\sigma)^p w \\
  &\leq \eta^{-p}\|M(f\sigma)\|_{L^p(w)}^p
  \leq \eta^{-p}\mathfrak{M}_p^p\|f\|_{L^p(\sigma)}^p.
\end{split}
\end{equation*}

We are left with
\begin{equation*}
\begin{split}
   \sum_{k\in\mathbb{Z}} 2^{kp} \sum_{Q\in\mathcal{Q}_k} w(E_k(Q))
   &\leq \sum_{k\in\mathbb{Z}} 2^{kp} \sum_{\substack{Q\in\mathcal{Q}_k \\ w(E_k(Q))<\delta w(Q)}} \delta w(Q) \\
    &\qquad+\sum_{k\in\mathbb{Z}} 2^{kp} \sum_{\substack{Q\in\mathcal{Q}_k \\ w(E_k(Q))\geq\delta w(Q)}} w(E_k(Q)),
\end{split}
\end{equation*}
and the first sum on the right is immediately dominated by
\begin{equation*}
   \delta \sum_{k\in\mathbb{Z}} 2^{kp} w(\mathbb{S}_{\natural}(f\sigma)>2^k)
   \lesssim \delta\|\mathbb{S}_{\natural}(f\sigma)\|_{L^p(w)}^p,
\end{equation*}
so this term can be absorbed after a suitable small choice of $\delta>0$ depending only on $p$.

Now consider one of the remaining sets $E_k(Q)$ for $Q\in\mathcal{Q}_k$ and $w(E_k(Q))\geq\delta w(Q)$. By the maximum principle, we can choose a linearization $\mathbb{L}=e^{i\vartheta(x)}\mathbb{S}_{\epsilon(x),\upsilon(x)}$ of $\mathbb{S}_{\natural}$ with $\upsilon(x)\leq\ell(Q)$, so that, for $x\in E_k(Q)$,
\begin{equation*}
\begin{split}
  4\cdot 2^k<\mathbb{S}_{\natural}(f\sigma)(x)
  &\leq 2\mathbb{L}(f\sigma)(x)+2^k+(\zeta+n+1)M(f\sigma)(x), \\
  &\leq 2\mathbb{L}(f\sigma)(x)+2^k+(\zeta+n+1)\eta 2^k \\
  &\leq 2\mathbb{L}(f\sigma)(x)+2\cdot 2^k,
\end{split}
\end{equation*}
by choosing 
\begin{equation*}
  \eta:=(\zeta+n+1)^{-1}=(2n+2)^{-1}\gtrsim\kappa^{-1}.
\end{equation*}
Hence
\begin{equation*}
  \mathbb{L}(f\sigma)(x)\geq 2^k\qquad\text{on }E_k(Q).
\end{equation*}
Notice that, by the disjointness of the sets $E_k(Q)\subset F_k(Q)$, we can globally define one linearization $\mathbb{L}$, which fulfills this condition on all $E_k(Q)$.

Thus, for $w(E_k(Q))\geq\delta w(Q)$, we have
\begin{equation*}
  2^{kp}
  \leq\Big(\frac{1}{w(E_k(Q))}\int_{E_k(Q)}\mathbb{L}(f\sigma)w\Big)^p
  \leq\delta^{-p}\Big(\frac{1}{w(Q)}\int_{E_k(Q)}\mathbb{L}(f\sigma)w\Big)^p,
\end{equation*}
where
\begin{equation*}
  \mathbb{L}=e^{i\vartheta(x)}\mathbb{S}_{\epsilon(x),\upsilon(x)},\qquad \upsilon(x)\leq\ell(Q)\text{ on }E_k(Q).
\end{equation*}
We have proven that (absorbing the term $\delta\|\mathbb{S}_{\natural}(f\sigma)\|_{L^p(\sigma)}^p$, and using $\eta^{-1}\lesssim\kappa$)
\begin{equation*}
  \|\mathbb{S}_{\natural}(f\sigma)\|_{L^p(\sigma)}^p
  \lesssim \kappa^p\mathfrak{M}_p^p\|f\|_{L^p(\sigma)}^p+\sum_{k\in\mathbb{Z}}\sum_{Q\in\mathcal{Q}_k}w(E_k(Q))
     \Big(\frac{1}{w(Q)}\int_{E_k(Q)}\mathbb{L}(f\sigma)w\Big)^p.
\end{equation*}
(The dependence on $\delta$ has been neglected, since this is in any case a function of $p$ only.)

\section{Strong-type estimates: the core}\label{s.nondoubling}

We are left to prove that
\begin{equation}\label{eq:mainTerm}
  \sum_k\sum_{Q\in\mathcal{Q}_k}w(E_k(Q))
    \Big|\frac{1}{w(Q)}\int_{E_k(Q)}\mathbb{L}(f\sigma)w\Big|^p
  \lesssim(\mathfrak{T}_p+\mathfrak{N}_p)^p\|f\|_{L^p(\sigma)}.
\end{equation}
We can make the additional assumption that all $k$ in this sum are of the same parity; after all, there are just two such sums.
By monotone convergence, we may also assume that all appearing cubes are contained in some maximal dyadic cube $\overline{Q}$. This allows to make the following construction:

\begin{definition}[Principal cubes]\label{def:prcubes}
We form the collection $\mathcal{G}$ of \emph{principal cubes} as follows. We let $\mathcal{G}_0:=\{\overline{Q}\}$ (the maximal cube that we consider), inductively
\begin{equation*}
  \mathcal{G}_k:=\bigcup_{G\in\mathcal{G}_{k-1}}\{G'\subset G:\mathbb{E}_{G'}^{\sigma}|f|>4\mathbb{E}_G^{\sigma}|f|,G'\text{ is a maximal such cube}\},
\end{equation*}
and then $\mathcal{G}:=\bigcup_{k=0}^{\infty}\mathcal{G}_k$. For any dyadic $Q$ ($\subset\overline{Q}$), we let
\begin{equation*}
  \Gamma(Q) :=\text{the minimal principal cube containing }Q.
\end{equation*}
\end{definition}

From the definition it follows that
\begin{equation*}
  \mathbb{E}_Q^{\sigma}|f|\leq 4\mathbb{E}_{\Gamma(Q)}^{\sigma}|f|.
\end{equation*}

We begin the analysis of \eqref{eq:mainTerm}. Recall that on $E_k(Q)$, we have $\mathbb{L}=\mathbb{L}(1_Q\,\cdot\,)$. Thus we may dualize and split the cube $ Q$ to the result that
\begin{equation}\label{eq:inAndOutOmega}
\begin{split}
  \int_{E_k(Q)}\mathbb{L}(f\sigma)w
  &=\int_Q\mathbb{L}^*(1_{E_k(Q)}w)f\sigma \\
  &=\int_{Q\setminus\Omega_{k+2}}\mathbb{L}^*(1_{E_k(Q)}w)f\sigma
    +\int_{Q\cap\Omega_{k+2}}\mathbb{L}^*(1_{E_k(Q)}w)f\sigma.
\end{split}
\end{equation}

\subsection{The part on $Q\setminus\Omega_{k+2}$}
The first term is easy to estimate:
\begin{equation*}
\begin{split}
  \Big|\int_{Q\setminus\Omega_{k+2}}\mathbb{L}^*(1_{E_k(Q)}w)f\sigma\Big|
  &\leq\|1_Q\mathbb{L}^*(1_{E_k(Q)}w)\|_{L^{p'}(\sigma)}\|1_{Q\setminus\Omega_{k+2}}f\|_{L^p(\sigma)} \\
  &\leq\mathfrak{T}_pw(Q)^{1/p'}\|1_{Q\setminus\Omega_{k+2}}f\|_{L^p(\sigma)},
\end{split}
\end{equation*}
and then
\begin{equation*}
\begin{split}
   \sum_k &\sum_{Q\in\mathcal{Q}_k}w(E_k(Q))\Big|\frac{1}{w(Q)}\int_{Q\setminus\Omega_{k+2}}\mathbb{L}^*(1_{E_k(Q)}w)f\sigma\Big|^p \\
   &\leq\mathfrak{T}_p^p \sum_k \sum_{Q\in\mathcal{Q}_k}w(E_k(Q))\frac{1}{w(Q)}\|1_{Q\setminus\Omega_{k+2}} f\|_{L^p(\sigma)}^p \\
   &\leq\mathfrak{T}_p^p \sum_k \sum_{Q\in\mathcal{Q}_k}\|1_{Q\setminus\Omega_{k+2}} f\|_{L^p(\sigma)}^p \\
   &\leq\mathfrak{T}_p^p \sum_k \|1_{\Omega_k\setminus\Omega_{k+2}} f\|_{L^p(\sigma)}^p 
     =\mathfrak{T}_p^p\|f\|_{L^p(\sigma)}^p,
\end{split}
\end{equation*}
recalling in the last step that the $k$ sum is over either odd or even $k$ only.

\subsection{The part on $Q\cap\Omega_{k+2}$}
We are left to estimate the integrals over $Q\cap\Omega_{k+2}$ as in the second term on the right of \eqref{eq:inAndOutOmega}. Using $Q\in\mathcal{Q}_k$ and $\Omega_{k+2}=\bigcup_{R\in\mathcal{Q}_{k+2}}R$, as well as the nestedness of the collections $\mathcal{Q}_k$, we have
\begin{equation*}
   \int_{Q\cap\Omega_{k+2}}\mathbb{L}^*(1_{E_k(Q)}w)f\sigma
   =\sum_{\substack{R\in\mathcal{Q}_{k+2}\\ R\subset Q}}\int_R \mathbb{L}^*(1_{E_k(Q)}w)f\sigma.
\end{equation*}
Now $E_k(Q)\subset\Omega_{k+2}^c\subset (R^{(\zeta)})^c$, where $\zeta=n+1$, so  $\mathbb{L}^*(1_{E_k(Q)}w)$ is in fact constant on $R\in\mathcal{Q}_{k+2}$ (by the `smoothness property' formulated in Lemma~\ref{lem:smoothness}); thus
\begin{equation*}
  \int_R \mathbb{L}^*(1_{E_k(Q)}w)f\sigma
  =\int_R \mathbb{L}^*(1_{E_k(Q)}w)\sigma\cdot \mathbb{E}_{R}^{\sigma}f.
\end{equation*}
Splitting into the cases according to the size of $\mathbb{E}_R^{\sigma}|f|$ relative to $\mathbb{E}_{\Gamma(Q)}^{\sigma}|f|$, we can thus estimate
\begin{equation}\label{eq:sizeOfAve}
\begin{split}
  \Big|\int_{Q\cap\Omega_{k+2}}\mathbb{L}^*(1_{E_k(Q)}w)f\sigma\Big| 
  &\leq 16\sum_{\substack{R\in\mathcal{Q}_{k+2}, R\subset Q\\ \mathbb{E}_R^{\sigma}|f|\leq16\mathbb{E}_{\Gamma(Q)}^{\sigma}|f|}}
       \int_R \big|\mathbb{L}^*(1_{E_k(Q)}w)\big|\sigma\cdot\mathbb{E}_{\Gamma(Q)}^{\sigma}|f| \\
  &\qquad+\sum_{\substack{R\in\mathcal{Q}_{k+2}, R\subset Q\\ \mathbb{E}_R^{\sigma}|f|>16\mathbb{E}_{\Gamma(Q)}^{\sigma}|f| }}
     \int_R \big|\mathbb{L}^*(1_{E_k(Q)}w)\big|\sigma\cdot\mathbb{E}_{R}^{\sigma}|f|.
\end{split}
\end{equation}

\subsection{The part with $\mathbb{E}_R^{\sigma}|f|\leq16\mathbb{E}_{\Gamma(Q)}^{\sigma}|f|$}
We estimate the first sum on the right of \eqref{eq:sizeOfAve}. As a first step, by the disjointness of the $R\in\mathcal{Q}_{k+2}$ we have
\begin{equation*}
   \sum_{\substack{R\in\mathcal{Q}_{k+2}, R\subset Q\\ \mathbb{E}_R^{\sigma}|f|\leq16\mathbb{E}_{\Gamma(Q)}^{\sigma}|f|}}
       \int_R \big|\mathbb{L}^*(1_{E_k(Q)}w)\big|\sigma
    \leq \int_Q\big|\mathbb{L}^*(1_{E_k(Q)}w)\big|\sigma.
\end{equation*}
Next, we make a manipulation involving random signs $\varepsilon_{k',Q'}$ indexed by pairs $k'\in\mathbb{Z}'$ and $Q'\in\mathcal{Q}_{k'}$. At almost every $x$,
\begin{align*}
\big|\mathbb{L}^*(1_{E_k(Q)}1_Qw)(x)\big| 
  &\leq\mathbb{E}_{\varepsilon}\Big|\sum_{k'}\sum_{Q'\in\mathcal{Q}_{k'}}\varepsilon_{k',Q'}
       \mathbb{L}^*(1_{E_{k'}(Q')}1_Qw)(x)\Big| \\
  &\leq\mathbb{E}_{\varepsilon}\Big|\mathbb{L}^*\Big(\sum_{k'}\sum_{Q'\in\mathcal{Q}_{k'}}\varepsilon_{k',Q'}
       1_{E_{k'}(Q')}1_Qw\Big)(x)\Big| \\
    &=:\mathbb{E}_{\varepsilon}\big|\mathbb{L}^*\big(\psi_{\varepsilon}1_Qw\big)(x)\big|
\end{align*}
where, by the disjointness of the sets $E_{k'}(Q')$, we have $\|\psi_{\varepsilon}\|_{\infty}\leq 1$ for any choice of the signs $\varepsilon_{k',Q'}$.

We can then compute
\begin{align*}
  \sum_k &\sum_{Q\in\mathcal{Q}_k}w(E_k(Q))\Big(\frac{1}{w(Q)}\int_Q\big|\mathbb{L}^*(1_{E_k(Q)}w)\big|\sigma
      \cdot\mathbb{E}_{\Gamma(Q)}^{\sigma}|f|\Big)^p \\
    &\leq\sum_{G\in\mathcal{G}}\big(\mathbb{E}_{G}^{\sigma}|f| \big)^p
       \sum_k\sum_{\substack{Q\in\mathcal{Q}_k\\ \Gamma(Q)=G}}w(E_k(Q))\Big(\frac{1}{w(Q)}\mathbb{E}_{\varepsilon}
           \int_Q\big|\mathbb{L}^*(\psi_{\varepsilon}1_Q w)\big|\sigma\Big)^p \\
           \noalign{\hfill (by the preceding considerations)}
   &\leq\sum_{G\in\mathcal{G}}\big(\mathbb{E}_{G}^{\sigma}|f| \big)^p\mathbb{E}_{\varepsilon}
       \sum_k\sum_{\substack{Q\in\mathcal{Q}_k\\ \Gamma(Q)=G}}w(E_k(Q))\Big(\frac{1}{w(Q)}
           \int_Q\big|\mathbb{L}^*(\psi_{\varepsilon}1_Q w)\big|\sigma\Big)^p \\
           \noalign{\hfill (by H\"older's inequality and the linearity of $\mathbb{E}_{\varepsilon}$)}
   &\leq\sum_{G\in\mathcal{G}}\big(\mathbb{E}_{G}^{\sigma}|f| \big)^p\mathbb{E}_{\varepsilon}
      \int_G\Big(\sup_{Q\subset G}\frac{1_Q}{w(Q)}
           \int_Q\big|\mathbb{L}^*(\psi_{\varepsilon}1_Q w)\big|\sigma\Big)^p w \\ 
           \noalign{\hfill (by the disjointness of the sets $E_k(Q$)}
   &\leq\sum_{G\in\mathcal{G}}\big(\mathbb{E}_{G}^{\sigma}|f| \big)^p\mathfrak{N}_p^p\sigma(G) \\
           \noalign{\hfill (by the non-standard testing condition)}
   &\lesssim\mathfrak{N}_p^p\|M_{\sigma}f\|_{L^p(\sigma)}^p\lesssim\mathfrak{N}_p^p\|f\|_{L^p(\sigma)}^p \\
   \noalign{\hfill (by definition of principal cubes and universal maximal inequality.)}
\end{align*}
 
\subsection{The part with $\mathbb{E}_R^{\sigma}|f|>16\mathbb{E}_{\Gamma(Q)}^{\sigma}|f|$}
We estimate the second term on the right of \eqref{eq:sizeOfAve}.
Inserting $\sigma(R)^{\pm 1/p}$ into this sum and using H\"older's inequality, we have
\begin{equation*}
\begin{split}
   & \sum_{\substack{R\in\mathcal{Q}_{k+2}, R\subset Q\\ \mathbb{E}_R^{\sigma}|f|>16\mathbb{E}_{\Gamma(Q)}^{\sigma}|f| }}
     \int_R \big|\mathbb{L}^*(1_{E_k(Q)}w)\big|\sigma\cdot\mathbb{E}_{R}^{\sigma}|f| \\
   &\leq\Big(\sum_{\substack{R\in\mathcal{Q}_{k+2}\\ R\subset Q}}
        \sigma(R)^{-p'/p}\Big[\int_R \big|\mathbb{L}^*(1_{E_k(Q)}w)\big|\sigma\Big]^{p'}\Big)^{1/p'}
        \Big(\sum_{\substack{R\in\mathcal{Q}_{k+2}, R\subset Q\\ \mathbb{E}_R^{\sigma}|f|>16\mathbb{E}_{\Gamma(Q)}^{\sigma}|f| }}
           \sigma(R)[\mathbb{E}_R^{\sigma}|f|]^p\Big)^{1/p},
\end{split}
\end{equation*}
where the first factor is further dominated by
\begin{equation*}
\begin{split}
   \Big(\sum_{\substack{R\in\mathcal{Q}_{k+2}\\ R\subset Q}}
        \int_R \big|\mathbb{L}^*(1_{E_k(Q)}w)\big|^{p'}\sigma\Big)^{1/p'}
   \leq\Big(\int_Q \big|\mathbb{L}^*(1_{E_k(Q)}w)\big|^{p'}\sigma\Big)^{1/p'}
   \leq \mathfrak{T}_pw(Q)^{1/p'}.
\end{split}
\end{equation*}
So altogether,  and appealing to the testing condition \eqref{e.tw2},
\begin{equation}\label{eq:finalPart}
\begin{split}
  \sum_k &\sum_{Q\in\mathcal{Q}_k}w(E_k(Q))
    \Big(\frac{1}{w(Q)}\sum_{\substack{R\in\mathcal{Q}_{k+2}, R\subset Q\\ \mathbb{E}_R^{\sigma}|f|>16\mathbb{E}_{\Gamma(Q)}^{\sigma}|f| }}
     \int_R \big|\mathbb{L}^*(1_{E_k(Q)}w)\big|\sigma\cdot\mathbb{E}_{R}^{\sigma}|f| \Big)^p \\
     &\leq\mathfrak{T}_p^p\sum_k \sum_{Q\in\mathcal{Q}_k}w(E_k(Q))\frac{1}{w(Q)}
        \sum_{\substack{R\in\mathcal{Q}_{k+2}, R\subset Q\\ \mathbb{E}_R^{\sigma}|f|>16\mathbb{E}_{\Gamma(Q)}^{\sigma}|f| }}
           \sigma(R)[\mathbb{E}_R^{\sigma}|f|]^p \\
     &\leq\mathfrak{T}_p^p\sum_k \sum_{Q\in\mathcal{Q}_k}
        \sum_{\substack{R\in\mathcal{Q}_{k+2}, R\subset Q\\ \mathbb{E}_R^{\sigma}|f|>16\mathbb{E}_{\Gamma(Q)}^{\sigma}|f| }}
           \sigma(R)[\mathbb{E}_R^{\sigma}|f|]^p
\end{split}
\end{equation}

The proof is completed by the following lemma:

\begin{lemma}
Any given cube $R$ appears at most once in the sum on the right of \eqref{eq:finalPart}. For any two cubes $R_1\subsetneq R_2$ appearing in this sum, we have $\mathbb{E}_{R_1}^{\sigma}|f|>4\mathbb{E}_{R_2}^{\sigma}|f|$.
\end{lemma}

\begin{proof}
We can prove both claims with one strike as follows. Let $R_1\subset R_2$ be cubes appearing in \eqref{eq:finalPart}, possibly equal. Thus for some cubes $Q_i$, we have
\begin{equation*}
  R_i\in\mathcal{Q}_{k_i+2},\qquad R_i\subset Q_i\in\mathcal{Q}_{k_i},\qquad \mathbb{E}_{R_i}^{\sigma}|f|>16\mathbb{E}_{\Gamma(Q_i)}^{\sigma}|f|.
\end{equation*}
 (In the equal case, we want to prove that $(k_1,Q_1)=(k_2,Q_2)$; in the unequal case, the estimate between the averages.)
Note that if $k_1=k_2$, then also $R_1=R_2$ and $Q_1=Q_2$, so there is nothing to prove. So let $k_1>k_2$, thus by the restriction on the $k$-sum, in fact $k_1\geq k_2+2$. Since the cubes $Q_1\in\mathcal{Q}_{k_1}$ and $R_2\in\mathcal{Q}_{k_2+2}$ intersect (on $R_1$) and $k_1\geq k_2+2$, the nestedness property implies that $Q_1\subset R_2$, and hence $\Gamma(Q_1)\subset \Gamma(R_2)$. Thus
\begin{equation*}
  \mathbb{E}_{R_1}^{\sigma}|f|>16\mathbb{E}_{\Gamma(Q_1)}^{\sigma}|f|
  \geq 16\mathbb{E}_{\Gamma(R_2)}^{\sigma}|f|\geq 4\mathbb{E}_{R_2}^{\sigma}|f|.
\end{equation*}
For $R_1=R_2$, this gives a contradiction, showing that the same $R$ cannot arise in the sum more than once. And for $R_1\subsetneq R_2$, this is precisely the asserted estimate.
\end{proof}

From the lemma it follows that
\begin{equation*}
\begin{split}
  \mathfrak{T}_p^p  \sum_k \sum_{Q\in\mathcal{Q}_k}&
        \sum_{\substack{R\in\mathcal{Q}_{k+2}, R\subset Q\\ \mathbb{E}_R^{\sigma}|f|>16\mathbb{E}_{\Gamma(Q)}^{\sigma}|f| }}
           \sigma(R)[\mathbb{E}_R^{\sigma}|f|]^p \\
      &\lesssim \mathfrak{T}_p^p\|M_{\sigma}f\|_{L^p(\sigma)}^p\lesssim\mathfrak{T}_p^p\|f\|_{L^p(\sigma)}^p,
\end{split}
\end{equation*}
and the proof is complete.

\section{Unweighted weak-type (1,1) inequalities} \label{s.w11}

We are now finished with the two-weight theory, and we start anew from a different corner of our proof scheme, see Figure~\ref{f.1}, the weak-type $(1,1)$ estimates for Haar shift operators. Again, we need bounds  that are effective in the complexity.   The following estimate was proved in \cite{1007.4330}*{Proposition 5.1}, with the additional observation concerning shifts with separated scales made in \cite{1010.0755}*{Theorem 5.2}.

\begin{definition}\label{d.scales}
We say that a shift $\mathbb{S}$ of complexity $\kappa$ has its \emph{scales separated} if all nonzero components $\mathbb{S}_Q$, $\mathbb{S}_{Q'}$ have $\log_2\ell(Q)\equiv\log_2\ell(Q')\mod(\kappa+1)$.  We likewise say that a subset $ \mathcal Q$ of the dyadic grid  
has \emph{scales separated} if $\log_2\ell(Q)\equiv\log_2\ell(Q')\mod(\kappa+1)$ for all $ Q,Q'\in \mathcal Q$.  
\end{definition}

\begin{proposition}\label{p.hytWeak}  
An  $ L ^2 $-bounded Haar shift operator $ \mathbb S $ of complexity $ \kappa $ maps $ L ^{1} (\mathbb R ^{d})$  into $L ^{1,\infty } (\mathbb R ^{d}) $ with norm at most  $C_d\kappa$. If $\mathbb{S}$ has scales separated, then the norm is at most $C_d$.
\end{proposition}


We need a strengthening of this Proposition for duals $\mathbb{L}^*$ of the all linearizations $\mathbb{L}$ of the maximal truncation $\mathbb{S}_{\natural}$ of $\mathbb{S}$. This type of result is not a classical one, and to prove it, we use a simplified version of the argument  of \cite{MR1783613} used to prove Carleson's Theorem on Fourier series \cite{MR0199631}. 

\begin{theorem}\label{t.weak} For an $ L ^2 $-bounded Haar shift operator $ \mathbb S $ of complexity $ \kappa $, we have
the following estimate,
uniform in $ \lambda >0$, and  compactly supported and bounded functions $ f $ on $ \mathbb R ^{d}$
\begin{align} 
\label{e.Lmax11}
\lambda \lvert  \{  \mathbb L ^{\ast}   f > \lambda \}\rvert  &\lesssim \kappa  \lVert f\rVert_{1} \,, 
\end{align}
where the  inequality holds uniformly in choice of the linearization $ \mathbb L$ 
of $\mathbb{S}_{\natural}$. If $\mathbb{S}$ has its scales separated, we have the complexity-free bound
\begin{equation}\label{eq:weak11separ}
  \lambda \lvert  \{  \mathbb L ^{\ast}   f > \lambda \}\rvert  \lesssim  \lVert f\rVert_{1} \,. 
\end{equation}
\end{theorem}
\subsection{Case $\int\mathbb{S}_Qf =0$}

We begin proving the weak type $(1,1)$ inequality with the additional hypothesis that the components of our shift $\mathbb{S}$ all satisfy $\int\mathbb{S}_Qf =0$. Note that this covers all cases relevant to us, except for the dual paraproduct case $\mathbb{S}_Q f=|Q|^{-1}\langle f, h_Q\rangle \mathbf 1_Q$, which will be taken up in the next subsection.

We start off with the \emph{Tree lemma}, with this terminology and notation adapted from \cite{MR1783613}. 

\begin{lemma}[Tree lemma]\label{treeLemma}
Suppose that $\mathcal{Q}$ is a collection of dyadic cubes, all contained in some $Q_0$. Then
\begin{equation*}
  \mathbb{L}_{\mathcal{Q}}f(x):=\sum_{\substack{Q\in\mathcal{Q}\\ \ell(Q)\leq\upsilon(x)}}\mathbb{S}_Qf(x)
\end{equation*}
satisfies
\begin{equation}\label{treeIneq}
  |\langle \mathbb{L}_{\mathcal{Q}}f, g\rangle|\lesssim\kappa\cdot\operatorname{size}_f(\mathcal{Q})
     \cdot\operatorname{dense}_g(\mathcal{Q})\cdot|Q_0|,
\end{equation}
where
\begin{equation*}
  \operatorname{size}_f(\mathcal{Q}) :=\sup_{Q\in\mathcal{Q}}\Big(\frac{1}{|Q|}\int_Q |f|^2\Big)^{1/2},\qquad
  \operatorname{dense}_g(\mathcal{Q}) :=\sup_{Q\in\mathcal{Q}}\sup_{Q'\supset Q}\frac{1}{|Q'|}\int_{Q'} |g|.
\end{equation*}
The factor $\kappa$ may be omitted in \eqref{treeIneq} if the scales of $\mathbb{S}$ are separated.
\end{lemma}

\begin{remark}\label{r.size} The notation of `size' and `density' are derived from \cite{MR1783613}. But size in this context 
is simpler, related to the maximal function $ (M f ^2) ^{1/2}  $.  The double supremum in the definition of density is 
essential for the inequalty \eqref{form1} below. 
\end{remark}
\begin{proof}
We may assume that $Q_0\in\mathcal{Q}$. For if not, let $Q_i$ be the maximal cubes in $\mathcal{Q}$, all contained in $Q_0$, and $\mathcal{Q}_i:=\{Q\in\mathcal{Q}:Q\subset Q_i\}$. Clearly the size and density of each $\mathcal{Q}_i$ is dominated by the corresponding number for $\mathcal{Q}$. Then we just write $\mathbb{L}_{\mathcal{Q}}=\sum_i\mathbb{L}_{\mathcal{Q}_i}$, use the estimate for each $\mathcal{Q}_i$, and sum up $\sum_i|Q_i|\leq|Q_0|$ in the end.

Also assume by approximation that the collection $\mathcal{Q}$ is finite. Let ${\mathcal P}$ consist of the minimal dyadic cubes $P \subset Q_0$ such that $P^{(1)}$ contains some element of ${\mathcal Q}$. The cubes in ${\mathcal P}$ then form a partition of $\bigcup \{Q : Q \in {\mathcal Q}\}$. Define the maximal operator ${\mathcal M}_{{\mathcal P}}$ by
\begin{displaymath}
  {\mathcal M}_{{\mathcal P}}\phi(x) := \sum_{P \in {\mathcal P}} |g(x)|\mathbf 1_{P \cap \{\ell(P^{(1)}) \leq \upsilon(x)\}} \sup_{Q \supset P} \mathbb E_{Q}|\phi|. 
\end{displaymath}
Note that if $P \in {\mathcal P}$, then $P^{(1)}$ contains a cube $Q \in {\mathcal Q}$ with $P \subset Q$, whence
\begin{displaymath} \mathbb E_{P}[|g|\mathbf 1_{\{\ell(P^{(1)}) \leq \upsilon(\cdot)\}}] \leq 2^{d}\mathbb E_{P^{(1)}}[|g|\mathbf 1_{\{\ell(P^{(1)}) \leq \upsilon(\cdot)\}}] \leq 2^{d}\textup{dense}_g({\mathcal Q}).
\end{displaymath}
From this we conclude a particular restricted strong-type inequality for $ \mathcal M _{\mathcal P}$:
\begin{align}
  \int_{Q_0} {\mathcal M}_{{\mathcal P}}\phi
  & \leq \sum_{P \in {\mathcal P}} |P| \cdot \mathbb E_{P}[|g|\mathbf 1_{\{\ell(P^{(1)}) \leq \upsilon(\cdot)\}}]\sup_{Q \supset P} \mathbb E_{Q}|\phi|\\
  & \leq 2^{d}\textup{dense}_g({\mathcal Q}) \sum_{P \in {\mathcal P}} |P| \sup_{Q \supset P} \mathbb E_{Q}|\phi| \\
  & \leq 2^{d}\textup{dense}_g({\mathcal Q}) \int_{Q_0} {\mathcal M}\phi\\
  &\label{form1} \leq 2^{d}\textup{dense}_g(\mathcal{Q}) |Q_0|^{1/2}\|{\mathcal M}\phi\|_{2} 
   \lesssim \textup{dense}_g({\mathcal Q}) |Q_0|^{1/2}\|\phi\|_{2}.
  \end{align}
 
Now we start to estimate the expression in \eqref{treeIneq}. For $x\in P\in\mathcal{P}$, we have
\begin{equation*}
  \mathbb{L}_{\mathcal{Q}}f(x)=\sum_{\substack{Q\in\mathcal{Q},\ Q\supsetneq P \\ \ell(Q)\leq\upsilon(x)}}\mathbb{S}_Qf(x)
  =1_{\{\ell(P)<\upsilon(x)\}}\sum_{\substack{Q\in\mathcal{Q},\ Q\supsetneq P \\ \ell(Q)\leq\upsilon(x)}}\mathbb{S}_Qf(x),
\end{equation*}
since the summation is empty unless $\ell(P)<\upsilon(x)$. Let $P_{\upsilon(x)}\supsetneq P$ be the unique dyadic cube with $\ell(P_{\upsilon(x)})=\upsilon(x)$. Then
\begin{equation*}
  1_{\{\ell(P)<\upsilon(x)\}}\sum_{\substack{Q\in\mathcal{Q},\ Q\supsetneq P \\ \ell(Q)\leq\upsilon(x)}}\mathbb{S}_Qf(x)
  =1_{\{\ell(P)<\upsilon(x)\}}\Big(\sum_{\substack{Q\in\mathcal{Q}\\ Q\supsetneq P}}\mathbb{S}_Qf(x)
      -\sum_{\substack{Q\in\mathcal{Q}\\ Q\supsetneq P_{\upsilon(x)}}}\mathbb{S}_Qf(x)\Big).
\end{equation*}
For any dyadic $R$ and $x\in R$, we have
\begin{equation}\label{eq:SQarrange}
   \sum_{\substack{Q\in\mathcal{Q}\\ Q\supsetneq R}}\mathbb{S}_Qf(x)
   =\mathbb{E}_R\Big(\sum_{Q\in\mathcal{Q}}\mathbb{S}_Qf\Big)
     +\sum_{\substack{Q\in\mathcal{Q}\\ R\subsetneq Q\subset R^{(\kappa)}}}\Big(\mathbb{S}_Qf(x)-\mathbb{E}_R(\mathbb{S}_Qf)\Big),
\end{equation}
which follows from the facts that $\mathbb{E}_R(\mathbb{S}_Qf)=0$ for $Q\subset R$, and $\mathbb{S}_Qf$ is constant on $R$ for $Q\supsetneq R^{(\kappa)}$. And here
\begin{equation*}
   \sum_{R\in\{P,P_{\upsilon(x)}\}}\Big|\mathbb{E}_R\Big(\sum_{Q\in\mathcal{Q}}\mathbb{S}_Qf\Big)\Big|
   \leq 2\sup_{R\supset P}\mathbb{E}_R|\mathbb{S}_{\mathcal{Q}}f|,\qquad
   \mathbb{S}_{\mathcal{Q}}f:=\sum_{Q\in\mathcal{Q}}\mathbb{S}_Q f.
\end{equation*}
For the second sum in \eqref{eq:SQarrange}, observe that we have 
\begin{equation*}
  |\mathbb{E}_R(\mathbb{S}_Qf)|\leq
  \|\mathbb{S}_Qf\|_{\infty}\leq\frac{1}{|Q|}\int_Q|f|\leq\Big(\frac{1}{|Q|}\int_Q|f|^2\Big)^{1/2}\leq\operatorname{size}_f(\mathcal{Q})
\end{equation*}
for each term, and there are at most $4\kappa$ terms altogether for both $R\in\{P,P_{\upsilon(x)}\}$.
If the scales of $\mathbb{S}$ are separated, then there is at most one nonzero $\mathbb{S}_Q$ with $R\subsetneq Q\subset R^{(\kappa)}$, so there are at most $4$ nonzero terms with the mentioned estimate, instead of $4\kappa$. This is the only place where the factor $\kappa$ enters into the argument, and the rest of the proof can be modified to the case of separated scales by simply substituting $1$ in place of $\kappa$. 

Substituting back, we have shown that
\begin{equation*}
  |\mathbb{L}_{\mathcal{Q}}f(x)|\leq \mathbf 1_{\{\ell(P)<\upsilon(x)\}}
    \big(2\sup_{R\supset P}\mathbb{E}_R|\mathbb{S}_{\mathcal{Q}}f|+4\kappa\operatorname{size}_f(\mathcal{Q})\big)\mathbf 1_{Q_0},\qquad \forall\,x\in P\in\mathcal{P},
\end{equation*}
and hence
\begin{equation*}
\begin{split}
  |\mathbb{L}_{\mathcal{Q}}f(x)\cdot g(x)|
  &\leq \sum_{P\in\mathcal{P}}|g(x)| \mathbf 1_{P \cap \{\ell(P^{(1)})\leq\upsilon(x)\}}\sup_{R\supset P}\mathbb{E}_R
      \big(2|\mathbb{S}_{\mathcal{Q}}f|+4\kappa\operatorname{size}_f(\mathcal{Q})\mathbf 1_{Q_0}\big) \\
  &=\mathcal{M}_{\mathcal{P}}\phi(x),\qquad \phi:=2|\mathbb{S}_{\mathcal{Q}}f|+4\kappa\operatorname{size}_f(\mathcal{Q})\mathbf 1_{Q_0}.
\end{split}
\end{equation*}
An application of \eqref{form1} then gives
\begin{equation*}
  |\langle \mathbb{L}_{\mathcal{Q}}f,g\rangle|\lesssim\operatorname{dense}_g(\mathcal{Q})\cdot|Q_0|^{1/2}\cdot\|\phi\|_2,
\end{equation*}
and
\begin{equation*}
\begin{split}
   \|\phi\|_2 &\lesssim\|\mathbb{S}_{\mathcal{Q}}(1_{Q_0}f)\|_2+\kappa\operatorname{size}_f(\mathcal{Q})\|\mathbf 1_{Q_0}\|_2 \\
     &\lesssim\|1_{Q_0}f\|_2+\kappa\operatorname{size}_f(\mathcal{Q})|Q_0|^{1/2}
     \leq (1+\kappa)\operatorname{size}_f(\mathcal{Q})|Q_0|^{1/2}.
\end{split}
\end{equation*}
This completes the proof, recalling that in the case of separated scales we can take $1$ in place of $\kappa$ above.
\end{proof}

 The next two Lemmas present decompositions of collections $ \mathcal Q$ relative to the two quantities of density and size.

\begin{lemma}\label{densityLemma}
Let ${\mathcal Q}$ be an arbitrary collection of dyadic cubes, $g\in L^1$, and $\delta>0$. Then $\mathcal{Q}=\mathcal{Q}'\cup\bigcup_j\mathcal{Q}_j$ where $\operatorname{dense}_g(\mathcal{Q}')\leq\delta\|g\|_1$ and $\mathcal{Q}_j=\{Q\in\mathcal{Q}:Q\subset Q_j\}$, where
\begin{equation*}
  \sum_j|Q_j|\leq\delta^{-1}.
\end{equation*}
\end{lemma}

\begin{proof}
Let $\mathcal{Q}' = \{Q \in \mathcal{Q}: \, \sup_{Q' \supset Q} |Q'|^{-1}\int_{Q'} |g| \le \delta\|g\|_1\}$. Then, by definition, there holds $\operatorname{dense}_g(\mathcal{Q}')\leq\delta\|g\|_1$.
Let $Q_j$ be the maximal elements of $\mathcal{Q} \setminus \mathcal{Q}'$, and $\mathcal{Q}_j$ be defined as in the statement. If $Q \in \mathcal{Q} \setminus \mathcal{Q}'$, then
$|R(Q)|^{-1} \int_{R(Q)} |g| > \delta\|g\|_1$ for some dyadic $R(Q) \supset Q$. Let $\mathcal{A}$ denote the maximal elements of $\{R(Q):\, Q \in \mathcal{Q} \setminus \mathcal{Q}'\}$.
We have
\begin{displaymath}
\sum_j |Q_j| = \sum_{Q \in \mathcal{A}} \sum_{j:\, R(Q_j) \subset Q} |Q_j| \le \sum_{Q \in \mathcal{A}} |Q| \le \sum_{Q \in \mathcal{A}} \frac{1}{\delta\|g\|_1} \int_Q |g| \le \frac{1}{\delta},
\end{displaymath}
where we used that $Q_j \subset R(Q_j)$, the disjointness of the cubes $Q_j$ and the disjointness of the cubes of $\mathcal{A}$.
\end{proof}

\begin{lemma}\label{sizeLemma}
Let ${\mathcal Q}$ be an arbitrary collection of dyadic cubes, $f\in L^2$, and $\sigma>0$. Then $\mathcal{Q}=\mathcal{Q}'\cup\bigcup_j\mathcal{Q}_j$ where $\operatorname{size}_f(\mathcal{Q}')\leq\sigma\|f\|_2$ and $\mathcal{Q}_j=\{Q\in\mathcal{Q}:Q\subset Q_j\}$, where
\begin{equation*}
  \sum_j|Q_j|\leq\sigma^{-2}.
\end{equation*}
\end{lemma}

\begin{proof}
Let $Q_j$ be the maximal elements of $\mathcal{Q}$ with $(|Q_j|^{-1}\int_{Q_j}|f|^2)^{1/2}>\sigma\|f\|_2$, if any, and $\mathcal{Q}_j$ be defined as in the statement. Then clearly $\mathcal{Q}':=\mathcal{Q}\setminus\bigcup_j\mathcal{Q}_j$ satisfies $\operatorname{size}_f(\mathcal{Q}')\leq\sigma\|f\|_2$, and
\begin{equation*}
  \sum_j|Q_j|\leq\sum_j\frac{1}{\sigma^2\|f\|_2^2}\int_{Q_j}|f|^2\leq\frac{1}{\sigma^2\|f\|_2^2}\|f\|_2^2=\frac{1}{\sigma^2},
\end{equation*}
since the $Q_j$ are disjoint by maximality.
\end{proof}

Inductive application of the previous two Lemmas leads to the following general result on decomposition of an arbitrary collection of cubes $ \mathcal Q$.

\begin{lemma}\label{decomposition}
Let ${\mathcal Q}$ be an arbitrary collection of dyadic cubes. Let $n_0$ be such that
\begin{align*}
\textup{dense}_g({\mathcal Q}) \le 2^{2n_0}\|g\|_1, \qquad
\textup{size}_f({\mathcal Q}) \le 2^{n_0}\|f\|_2.
\end{align*}
We can write a disjoint union
\begin{displaymath}
  {\mathcal Q} = \mathcal{Q}^{-\infty}\cup\bigcup_{k=-\infty}^{n_0}\bigcup_j {\mathcal Q}^k_j
\end{displaymath}
where
{\renewcommand{\theenumi}{\roman{enumi}}
\begin{enumerate}
\item $\textup{dense}_g({\mathcal Q}^k_j) \le 2^{2k}\|g\|_1$,
\item $\textup{size}_f({\mathcal Q}^k_j) \le 2^k\|f\|_2$,
\item all cubes in $\mathcal{Q}^k_j$ are contained in one $Q^k_j$, and
\begin{equation*}
  \sum_j|Q^k_j|\leq 8\cdot 2^{-2k},
\end{equation*}
\item both $g$ and $f$ vanish almost everywhere on all $Q\in\mathcal{Q}^{-\infty}$.
\end{enumerate}
}
\end{lemma}

\begin{proof}
Using the previous Lemmas, we first split $\mathcal{Q}=\mathcal{Q}'\cup\bigcup_j\mathcal{Q}'_j$ where $\operatorname{dense}_g(\mathcal{Q}')\leq 2^{2(n_0-1)}\|g\|_1$, and the cubes of $\mathcal{Q}'_j$ are contained in $Q'_j$ with $\sum_j|Q_j'|\leq 2^{-2(n_0-1)}$. Next, $\mathcal{Q}'=\mathcal{Q}''\cup\bigcup_j\mathcal{Q}_j''$ where $\operatorname{size}_f(\mathcal{Q}'')\leq 2^{n_0-1}\|f\|_2$, and the cubes of $\mathcal{Q}_j''$ are contained in $Q_j''$ with $\sum_j|Q_j''|\leq 2^{-2(n_0-1)}$. We re-enumerate the collections $\mathcal{Q}_j'$ and $\mathcal{Q}_j''$ as $\mathcal{Q}_j^{n_0}$, similarly for the containing cubes which satisfy $\sum_j|Q_j^{n_0}|\leq 8\cdot 2^{-2n_0}$. Since $\mathcal{Q}^{n_0}_j\subset\mathcal{Q}$, these have density and size as required, and $\mathcal{Q}''\subset\mathcal{Q}'$ has $\operatorname{dense}_g(\mathcal{Q}'')\leq 2^{2(n_0-1)}$ and $\operatorname{size}_f(\mathcal{Q}'')\leq 2^{n_0-1}$ by construction. We may thus iterate with $(\mathcal{Q},n_0)$ replaced by $(\mathcal{Q}'',n_0-1)$. If some cube $Q\in\mathcal{Q}$ is not chosen to any $\mathcal{Q}^k_j$ at any phase of the iteration, this means that both $\int_Q|g|=\int_Q|f|^2=0$; these cubes constitute the collection $\mathcal{Q}^{-\infty}$.
\end{proof}

Now we are prepared to prove the weak type $(1,1)$ inequality for $\mathbb L^{\ast}$. 

\begin{proof}[Proof of \eqref{e.Lmax11}.]
We may assume that $\lambda = 1$ and that $|E| < \infty$, where $E := \{|\mathbb L^{\ast}g| > 1\}$. 
And consider the set  $G = E \cap \{Mg \le C_d\|g\|_1/|E|\}$. Fixing $C_d$ large enough (depending on the dimension $d$ only), there holds 
\begin{displaymath}
  \frac12|E| \leq |G|
  \le\int_G|\mathbb{L}^*g|=\langle f,\mathbb{L}^*g\rangle
  =\langle \mathbb{L}_{\mathcal{Q}}f, g\rangle,
\end{displaymath}
where $|f|=\mathbf 1_G$, and $\mathcal{Q} = \{Q \in {\mathcal D}: \, Q \cap G \ne \emptyset\}$; note that $\mathbb{S}_Q f=0$ unless $Q\in\mathcal{Q}$.

{The point of passing to the supplementary set $ G$ is that we have the estimate } 
\begin{displaymath}
\textup{dense}_g({\mathcal Q}) \le \frac{C_d}{|E|}\|g\|_1.
\end{displaymath}
{And, as $|f|=1_G$, it also follows that $\textup{size} _{f} ({\mathcal Q}) \leq 1=  |G|^{-1/2}\|f\|_2$.} 
 We apply Lemma \ref{decomposition} to ${\mathcal Q}$, which yields the corresponding decomposition of $\mathbb{L}_{\mathcal{Q}}$. Note that $\langle \mathbb{L}_{\mathcal{Q}^{-\infty}}f,g\rangle=0$. Hence
 \begin{equation*}
 \begin{split}
  |\{|\mathbb{L}^*g|>1\}|
  &=|E|\lesssim |\langle\mathbb{L}_{\mathcal{Q}}f,g\rangle|
    \leq\sum_{k,j}|\langle\mathbb{L}_{\mathcal{Q}^k_j}f,g\rangle| \\
  &\lesssim\kappa\sum_{k,j}\operatorname{size}_j(\mathcal{Q}^k_j)\cdot\operatorname{dense}_g(\mathcal{Q}^k_j)\cdot|Q^k_j| \\
  &\lesssim\kappa\sum_k 2^k\|f\|_2\cdot\min\{2^{2k},|E|^{-1}\}\|g\|_1\cdot\sum_j|Q^k_j| \\
  &\lesssim\kappa\sum_k 2^k|E|^{1/2}\cdot\min\{2^{2k},|E|^{-1}\}\|g\|_1\cdot 2^{-2k}
  \lesssim\kappa\|g\|_1.
\end{split}
\end{equation*}
If the scales of $\mathbb{S}$ are separated, the factor $\kappa$ does not appear in the application of the Tree Lemma~\ref{treeLemma}, and we obtain instead that $|\{|\mathbb{L}^*g|>1\}|\lesssim\|g\|_1$, as claimed. This completes the proof.
\end{proof}

\subsection{The case $\int\mathbb{S}_Qf\neq 0$}

As mentioned, this only appears in the dual paraproduct case where $\mathbb{S}_Qf=|Q|^{-1}\langle f,h_Q\rangle\mathbf 1_Q$, where $h_Q$ is a Haar function on $Q$. In this case, the operator $\mathbb{L}$ has the form
\begin{equation*}
  \mathbb{L} ^{\ast} g(x)=\sum_{Q\in\mathcal{D}}|Q|^{-1}h_Q(x)\langle\mathbf 1_{Q\cap\{\ell(Q)\leq\upsilon(Q)\}},g\rangle
  =\sum_{Q\in\mathcal{D}}|Q|^{-1}h_Q(x)a_Q\langle\mathbf 1_Q,|g|\rangle
  =:\tilde{\mathbb{S}}|g|(x),
\end{equation*}
where $|a_Q|\leq 1$ are some numbers, of course depending on $g$. However, the new shifts $\tilde{\mathbb{S}}f(x):=\sum_{Q\in\mathcal{D}}|Q|^{-1}h_Q(x) a_Q\langle\mathbf 1_Q,f\rangle$ are uniformly bounded on $L^2$ with respect to the choice of the $a_Q$, hence by the weak-type estimate for untruncated shifts, also uniformly bounded from $L^1$ to $L^{1,\infty}$. Thus
\begin{equation*}
  \|\mathbb{L}^*g\|_{L^{1,\infty}}=\|\tilde{\mathbb{S}}g\|_{L^{1,\infty}}\lesssim\|g\|_{L^1},
\end{equation*}
and we are done in this case as well.

\section{Distributional Estimates} \label{s.dist}

\subsection{A John-Nirenberg Estimate}

We recall this formulation from \cite{1010.0755}*{Lemma 5.5}.  
 Let  $\mathcal D_\kappa$ be  a scales seperated grid, as defined in Definition~\ref{d.scales}.

\begin{definition}
\label{df:maxf}
Let $\{ \phi_ Q \;:\; Q\in \mathcal D_\kappa \}$ be a collection of functions such that $\phi_ Q$ is supported on $Q$ and is constant on its $ \mathcal D_\kappa$-children. For $R_0\in\mathcal D_\kappa$ let   $\phi^*_{R_0}$ be a maximal function
\begin{equation*}
\phi^*_{R_0} (x) := \sup_{Q\in\mathcal D_\kappa: Q \ni x} \ \Bigl| \sum_{\substack{R\in\mathcal D_\kappa\\  Q\subset R \subset R_0}} \phi_ R(x) \Bigr|. 
\end{equation*}
\end{definition}

\begin{lemma}
\label{l.JohnNir}
Let $ \{ \phi_ Q \;:\; Q\in\mathcal D_\kappa \} $ be a collection of functions such that 
\begin{enumerate}
    \item $\phi_ Q$ is supported on $Q$ and constant on the  $ \mathcal D _{\kappa }$-children of $Q$;
    \item $\|\phi_ Q \|_\infty\le 1$; 
    \item There exists $\delta\in(0,1)$ such that for all cubes $R\in \mathcal D_\kappa$ 
\begin{equation*}
    \bigl| \Bigl\{  x\in R :   \phi^*_ R(x)  > 1 \Bigr\} \bigr| \le \delta |R| \,.
\end{equation*}
    \end{enumerate}
Then for all $R\in\mathcal D_\kappa$ and for all $t>1$
\begin{equation*}
\bigl| \Bigl\{  x\in R :   \phi^*_ R(x)  > t \Bigr\} \bigr| \le
\delta^{(t-1)/2} |R|\,.
\end{equation*}
\end{lemma}

\subsection{The Corona And Distributional Estimates}

We need the important definition of the stopping cubes, and Corona Decomposition.    
The grid $ \mathcal D _{\kappa } \subset \mathcal D$ has scales separated, as in Definition~\ref{d.scales}. 

\begin{definition}\label{d.stopping} Let $ w $ be a weight and $ Q \in \mathcal D _{\kappa }$.  We set the 
\emph{$ w$-stopping children} $ \mathcal S (Q)$ to be the maximal subcubes $ Q' \subset Q$ with $ Q' \in \mathcal D _{\kappa }$, so that 
\begin{equation}\label{e.4}
\frac {w (Q')} {\lvert  Q'\rvert } \ge 4 \frac {w (Q)} {\lvert  Q\rvert }. 
\end{equation}
Setting $ \mathcal S_0 = \mathcal S (Q)$, and inductively setting $ \mathcal S _{j+1} := \bigcup _{Q' \in \mathcal S _{j}} \mathcal S (Q')$, 
we refer to $ \mathcal S := \bigcup _{j=0} ^{\infty } \mathcal S _{j}$ as the \emph{$w$-stopping cubes of $ Q$.} 

The \emph{$ w$-Corona Decomposition} of a collection of cubes $ \mathcal Q$ with $ Q' \subset Q$ for all $ Q' \in \mathcal Q$, 
consists of the $ w$-stopping cubes $ \mathcal S$, and collections of cubes $ \{ \mathcal P (S) \;:\; S\in \mathcal S \}$
 so that  (1) the collections $  \mathcal P (S)$ form a disjoint decomposition of $ \mathcal Q$, and  (2) for all $ S\in \mathcal S$, 
 and $ Q' \in \mathcal P (S)$, we have that $ S$ is the minimal element of $ \mathcal S$ that contains $ Q'$.  In particular, we have 
 \begin{equation} \label{e.4ps}
\frac {w (Q')} {\lvert  Q'\rvert } < 4 \frac {w (S)} {\lvert  S\rvert } \,, \qquad Q'\in \mathcal P (S)\,, \ S\in \mathcal S \,. 
\end{equation}

\end{definition}

The previous definitions make sense for any weight $ w$.  Specializing to $ w\in A_p$ leads to the following elementary, and essential, observations.  The first is a familiar inequality, showing that an $ A_p $ weight cannot be too concentrated.  

\begin{lemma}\label{l.Ainfty}  Let $ w \in A_p$, $ Q$ is a cube and $ E \subset Q$.  We then have 
\begin{equation}\label{e.Ainfty}
\Biggl[\frac {\lvert  E\rvert } {\lvert  Q\rvert } \Biggr] ^{p} \lVert w\rVert_{A_p} ^{-1} \le \frac {w (E)} {w (Q)} \,. 
\end{equation}
\end{lemma}

\begin{proof} 
The property that $ w >0$ a.e. allows us to write 
\begin{align*}
\frac {\lvert E\rvert} {\lvert Q\rvert}  & = 
\frac {\int _{E} w ^{1/p} (x) w (x) ^{-1/p} \; dx } {\lvert  Q\rvert } 
\\
& \le  \frac {w (E) ^{1/p}  \sigma (Q) ^{1/p'}} {\lvert  Q\rvert } 
\\
& = \Bigl[  \frac {w (E) } {w (Q)} \Bigr] ^{1/p}   \frac {w (Q) ^{1/p} \sigma (Q) ^{1/p'}} {\lvert  Q\rvert } 
\end{align*}
which proves the Lemma. 
\end{proof}

The second is a direct application of the previous assertion.

\begin{lemma}\label{l.stoppingSum} Let $ w \in A_p$, $ Q$ a cube, and $ \mathcal S$ the $ w$-stopping cubes for $ Q$.  Then, we have 
\begin{equation}\label{e.stoppingSum}
\sum_{S\in \mathcal S} w (S) \lesssim \lVert w\rVert_{A_p} w (Q) \,. 
\end{equation}
\end{lemma}

\begin{proof}
It follows from \eqref{e.4} that we have that the union of the stopping children $ \mathcal S (Q)$ has Lebesgue measure at 
most $ \tfrac 14 \lvert  Q\rvert $.  Applying \eqref{e.Ainfty} to the set $ E = Q \backslash \bigcup _{S\in \mathcal S (Q)} S$, it 
follows that 
\begin{align*}
\frac {w (E)} {w (Q)} & \ge (3/4) ^{p} \lVert w\rVert_{A_p} ^{-1} \,.  
\end{align*}
Thus, we have 
\begin{equation*}
\sum _{S\in \mathcal S (Q)} w (S) \le (1 - c \lVert w\rVert_{A_p} ^{-1}) w (Q) \,. 
\end{equation*}
Recalling the inductive definition of the $ w$-stopping children, we see that our estimate holds.  
\end{proof}

\subsection{The Distributional Estimates}

We combine the John-Nirenberg and Corona Decomposition to obtain the crucial distributional estimates: The operator
 $ \mathbb L$, decomposed using the Corona, satisfies exponential distributional inequalities.  We will then strongly use the 
exponential moments to control certain $ L ^{p}$ norms in the next section.  
Our estimates are two fold.   The first  inequalities involve the Lebesgue measure, 
which are the important intermediate step to obtain the second inequalities for the $ \sigma $ measure.

\begin{definition}\label{d.adapted} For $ w \in A_p$,  and integers $ \alpha , \kappa \in \mathbb N  $, we will say that 
a collection $ \mathcal Q $ of cubes is $ (w, \alpha , \kappa )$-adapted if these conditions hold.  
First, for any $ Q, Q'\in \mathcal Q$, we have $ \log_2\ell (Q) = \log_2 \ell (Q') \mod (\kappa + 1)$. Second, we have 
\begin{equation} \label{e.ApFixed}
2 ^{\alpha } \le \frac {w (Q)} {\lvert  Q\rvert } 
\Bigl[  \frac {\sigma  (Q)} {\lvert  Q\rvert }  \Bigl] ^{p - 1}  < 2 ^{\alpha +1} \,, \qquad Q\in \mathcal Q \,. 
\end{equation}
We only need to consider $ 0 \le \alpha \le \lceil \log_2 \lVert w\rVert_{A_p} \rceil$.  
\end{definition}

\begin{lemma}\label{l.distribL}
 Let $ w \in A_p$, with dual measure $ \sigma $,  and let $ Q $ be a cube.  
For integers $ \alpha ,\kappa $, let $ \mathcal Q$ be a collection of cubes contained in $ Q$ which are $(\sigma ,\alpha ,\kappa  )$-adapted.  
Construct the $ w  $-Corona Decompositions $ \mathcal P (S)$, $ S\in \mathcal S$ of 
$ \mathcal Q$.   Let $ \mathbb S$ be a Haar Shift operator of complexity $ \kappa $, with $ \lVert \mathbb S \rVert_{ L ^2 \mapsto L ^2 } =1 $.  We have these estimates, uniform in  $ t\ge 1$, choices of linearizations $ \mathbb L$, 
functions $ \varphi $ with $ \lVert \varphi \rVert_{\infty } \le 1$, 
and $ S\in \mathcal S$: 
\begin{align}
\label{e.Lexp1} \bigl\lvert \bigl\{ x \in S \;:\; \mathbb L^{\ast} _{\mathcal P (S)} ( \mathbf 1_{S} \varphi w ) 
\ge \mathsf K t \frac {w (S)} {\lvert  S\rvert }
\bigr\}\bigr\rvert & \lesssim 2^{-t} \lvert S\rvert \,,  \\ 
\label{e.Lexp2} \sigma\bigl(\bigl\{ x \in S \;:\; \mathbb L^{\ast} _{\mathcal P (S)} ( \mathbf 1_{S} \varphi w ) 
\ge \mathsf K t \frac {w (S)} {\lvert  S\rvert }
\bigr\} \bigr)  &\lesssim 2^{-t} \sigma (S) \,. 
\end{align}
where $\mathsf K$ is a dimensional constant.
\end{lemma}

Observe that the condition that $\mathcal Q$ be $(\sigma,\alpha,\kappa)$-adapted implies in particular that all $\mathbb{L}_{\mathcal{Q}'}$, with $\mathcal{Q}'\subset\mathcal{Q}$, can be viewed as linearizations of shifts with separated scales; this will place the stronger conclusion \eqref{eq:weak11separ} of Theorem~\ref{t.weak} at our disposal, allowing us to get the stated complexity-free estimates \eqref{e.Lexp1} and \eqref{e.Lexp2}. The proof given below follows almost to the letter the proof of \cite{1010.0755}*{Lemma 5.6}.



\subsection{Proof of the distributional estimate for the Lebesgue measure}
We begin with the bound \eqref{e.Lexp1}. We aim to apply the John-Nirenberg estimate, Lemma \ref{l.JohnNir}. To this end, write
\begin{align*}
  \mathbb{L}^{\ast}_{\mathcal{P}(S)}(\mathbf{1}_{S}\varphi w)(y)
  & = \sum_{Q \in \mathcal{P}(S)} \int_{Q \cap \{\varepsilon(x) \leq |Q| \leq \upsilon(x)\}} e^{i\vartheta(x)}s_{Q}(x,y)\mathbf{1}_{S}(x)\varphi(x)w(x) \, dx\\
  & =: \sum_{Q \in \mathcal{P}(S)} \phi_{Q}(y),
 \end{align*} 
and 
\begin{displaymath} \phi_{R}^{\ast}(x) := \sup_{Q \in \mathcal{D}_{\kappa} \colon Q \ni x} \left|\sum_{Q' \in \mathcal{P}(S) \colon Q \subsetneq Q' \subset R} \phi_{Q'}(x)\right|.
 \end{displaymath}
for $R \in \mathcal{D}_{\kappa}$ (note that $\phi_{R}^{\ast} \equiv 0$, if $R \not\subset S$). Then $|\mathbb{L}^{\ast}_{\mathcal{P}(S)}(\mathbf{1}_{S}\varphi w)| \leq \phi^{\ast}_{S}$, so that it suffices to prove \eqref{e.Lexp1} for $\phi^{\ast}_{S}$ instead of $\mathbb{L}^{\ast}_{\mathcal{P}(S)}(\mathbf{1}_{S}\varphi w)$. The condition $(1)$ of Lemma \ref{l.JohnNir} is clear: the function $\phi_{Q}$ is supported on $Q$, since the kernel $s_{Q}$ is supported on $Q \times Q$. Also, $\phi_{Q}$ is constant on $\mathcal{D}_{\kappa}$-children of $Q$. The condition $(2)$, that $\|\phi_{Q}\|_{\infty} \leq 1$, need not hold for cubes $Q \in \mathcal{P}(S)$, so we have to split these cubes into countably many subfamilies. For $\beta \in \mathbb{Z}_{+} = \{0,1,2,\ldots\}$, let $\mathcal{P}_{\beta}(S)$ consist of all cubes $Q \in \mathcal{P}(S)$ such that
\begin{displaymath} 4^{-\beta}\frac{w(S)}{|S|} \leq \frac{w(Q)}{|Q|} < 4^{-\beta + 1}\frac{w(S)}{|S|}.
 \end{displaymath}
For every $Q \in \mathcal{P}(S)$ it holds that
\begin{displaymath} \frac{w(Q)}{|Q|} < 4\frac{w(S)}{|S|},
 \end{displaymath} 
which means that every $Q \in \mathcal{P}(S)$ is contained in $\mathcal{P}_{\beta}(S)$ for some $\beta \in \mathbb{Z}_{+}$. Defining $\phi_{R,\beta}^{\ast}$ in the same manner as $\phi_{R}^{\ast}$, only replacing $\mathcal{P}(S)$ by $\mathcal{P}_{\beta}(S)$ in the definition, we also have
\begin{displaymath} \phi_{S}^{\ast} \leq \sum_{\beta \in \mathbb{Z}_{+}} \phi_{S,\beta}^{\ast}.
 \end{displaymath}
This will allow us to prove \eqref{e.Lexp1} for the functions $\phi_{S,\beta}^{\ast}$ separately. If $Q \in \mathcal{P}_{\beta}(S)$, we have
\begin{displaymath} |\phi_{Q}(y)| \leq \frac{1}{|Q|} \int_{Q \cap \{\varepsilon(x) \leq |Q| \leq \upsilon(x)\}} w(x) \, dx \leq 4^{-\beta + 1}\frac{w(S)}{|S|},
\end{displaymath}
by definition of $\mathcal{P}_{\beta}(S)$. Thus condition $(2)$ of Lemma \ref{l.JohnNir} holds for the normalized functions $4^{\beta - 1}|S|w(S)^{-1}\phi_{Q}$, $Q \in \mathcal{P}_{\beta}(S)$. In order to verify condition $(3)$ for $\phi_{R,\beta}^{\ast}$, we need the weak-type (1,1) inequality for $\mathbb{L}^{\ast}$ obtained above. To use this inequality, we first have to show that the set $\{x \in R : \phi_{R,\beta}^{\ast}(x) > \lambda\}$ is a subset of
\begin{displaymath} \{x \in R : |\mathbb{L}^{\ast}_{\mathcal{Q}_{\beta}(R)}(\mathbf{1}_{R}\varphi w)(x)| > \lambda\} 
\end{displaymath}
for some subcollection $\mathcal{Q}_{\beta}(R) \subset \mathcal{P}_{\beta}(S)$. Let $Q_{1},Q_{2},\ldots$ be the maximal subcubes $Q \subset R$ in $\mathcal{D}_{k}$ satisfying
\begin{displaymath} \left|\sum_{Q' \in \mathcal{P}_{\beta}(S) \colon Q \subsetneq Q' \subset R} \phi_{Q'}(x) \right| > \lambda.
\end{displaymath}
Note that every point in $\{x \in R : \phi_{R,\beta}^{\ast}(x) > \lambda\}$ is contained in exactly one $Q_{j}$ and, in fact, $\{x \in R : \phi_{R,\beta}(x) > \lambda\} = \bigcup_{j} Q_{j}$. Then let 
\begin{displaymath} \mathcal{Q}_{\beta}(R) := \bigcup_{j = 1}^{\infty} \{Q \in \mathcal{P}_{\beta}(S) : Q_{j} \subsetneq Q \subset R\}.
\end{displaymath}
Now suppose $\phi_{R,\beta}^{\ast}(x) > \lambda$, and let $Q_{j}$ be the unique cube containing $x$. Then, by definition,
\begin{align*}
  \lambda
  & < \left|\sum_{Q \in \mathcal{P}_{\beta}(S) \colon Q_{j} \subsetneq Q \subset R} \phi_{Q}(x) \right|\\
  & =: \left|\sum_{Q \in \mathcal{P}_{\beta}(S) \colon Q_{j} \subsetneq Q \subset R}
          \int_{Q \cap \{\varepsilon(y) \leq |Q| \leq \upsilon(y)\}} e^{i\vartheta(y)}s_{Q}(y,x)\mathbf{1}_{S}(y)\varphi(y)w(y) \, dy \right|.
\end{align*}
Since every cube in the sum is a subset of $R$, we may replace $\mathbf{1}_{S}$ by $\mathbf{1}_{R}$ above. The result is precisely $|\mathbb{L}_{\mathcal{Q}_{\beta}(R)}^{\ast}(\mathbf{1}_{R}\varphi w)(x)|$.

Let $\mathsf K\geq 1$ be the dimensional constant of Theorem~\ref{t.weak} for the weak-type inequality \eqref{eq:weak11separ} for shifts with scales separated. Recalling that this separation is satisfied in our situation, we have in particular that 
\begin{displaymath}
  |\{|\mathbb{L}^{\ast}_{\mathcal{Q'}}f| > \lambda\}| \leq \frac{\mathsf{K}\|f\|_{1}}{\lambda}, \qquad \lambda > 0,
\end{displaymath}
for any collection $\mathcal{Q}'\subset\mathcal{Q}$.
Choosing $\lambda = 2 \cdot |S|^{-1}w(S)4^{-\beta + 1}\mathsf{K}$ now yields
\begin{align*} |\{x \in R : \phi_{R,\beta}^{\ast}(x) > \lambda\}| & \leq |\{x \in R : |\mathbb{L}^{\ast}_{\mathcal{Q}_{\beta}(R)}(\mathbf{1}_{R}\varphi w)(x)| > \lambda\}|\\
& \leq \frac{\mathsf{K}|S|\|\mathbf{1}_{R}w\|_{1}}{2 \cdot w(S)4^{-\beta + 1}\mathsf{K}} = \frac{|S|\cdot w(R)}{2\cdot w(S)4^{-\beta + 1}}. 
\end{align*} 
If $R \in \mathcal{P}_{\beta}(S)$, we immediately obtain $|\{x \in R : \phi_{R,\beta}^{\ast}(x) > \lambda\}| \leq |R|/2$. However, Lemma \ref{l.JohnNir} requires the same estimate for all $R \in \mathcal{D}_{\kappa}$. Let $R_{1},R_{2},\ldots$ be the maximal cubes of $\mathcal{P}_{\beta}(S)$ inside $R$. Then $\phi_{R,\beta}^{\ast}$ is supported on the disjoint union of the cubes $R_{j}$, and $\phi_{R,\beta}^{\ast}(x) = \phi_{R_{j},\beta}^{\ast}(x)$ for $x \in R_{j}$. Hence
\begin{align*} |\{x \in R : \lambda^{-1}\phi_{R,\beta}^{\ast}(x) > 1\}| & = \sum_{j = 1}^{\infty} |\{x \in R_{j} : \phi_{R_{j},\beta}^{\ast}(x) > \lambda\}| \leq \sum_{j = 1}^{\infty} \frac{|R_{j}|}{2} \leq \frac{|R|}{2}. 
\end{align*}
Since $\lambda \geq |S|^{-1}w(S)4^{-\beta + 1}$, we still have $\|\lambda^{-1}\phi_{Q}\|_{\infty} \leq 1$, so that finally the functions $\lambda^{-1}\phi_{Q}$, and the related $\lambda^{-1}\phi_{R,\beta}^{\ast}$, satisfy all requirements of Lemma \ref{l.JohnNir}. We conclude that
\begin{displaymath} |\{x \in S : \phi_{S,\beta}^{\ast}(x) > t\lambda\}| \leq 2^{-(t - 1)/2}|S|, \qquad t > 1. 
\end{displaymath} 
Note that the inequality is trivial for $0 < t < 1$, so that it holds for all $t > 0$. Recalling the definition of $\lambda$ and rescaling $t$ we can rewrite the previous as
\begin{equation}\label{de:form2} \left|\left\{x \in S : \phi_{S,\beta}^{\ast}(x) > 16t\frac{w(S)}{|S|}\right\}\right| \leq \sqrt{2} \cdot 2^{-t4^{\beta}/\mathsf{K}}|S|, \qquad t > 0. 
\end{equation}
To finish the proof of \eqref{e.Lexp1}, we use $\phi_{S}^{\ast} \leq \sum \phi_{S,\beta}^{\ast}$ to estimate
\begin{align*} |S|^{-1} & \left|\left\{x \in S : \phi_{S}^{\ast}(x) > 16t\frac{w(S)}{|S|}\right\}\right|\\
& \leq |S|^{-1} \sum_{\beta = 0}^{\infty} \left|\left\{x \in S : \phi_{S,\beta}^{\ast}(x) > 16 \cdot 2^{-\beta - 1}t\frac{w(S)}{|S|}\right\}\right|\\
& \leq \sum_{\beta = 0}^{\infty} \sqrt{2} \cdot 2^{-t \cdot 2^{\beta - 1}/\mathsf{K}} \stackrel{(\ast)}{\leq} \sqrt{2}\sum_{\beta = 0}^{\infty} 2^{-t/2\mathsf{K} - \beta} = 2\sqrt{2} \cdot 2^{-t/2\mathsf{K}}. 
\end{align*} 
The inequality $(\ast)$ only holds in case $t \geq 2\mathsf{K}$, but for $t < 2\mathsf{K}$ the inequality just obtained is trivial. Replacing $t$ by $\mathsf{K}t$ proves \eqref{e.Lexp1}.


\subsection{Proof of the distributional estimate for the $\sigma$ measure}
In order to prove \eqref{e.Lexp2}, we need the assumption that our cube family is $(\sigma,\alpha,\kappa)$-adapted. Namely, recalling also the definition of $\mathcal{P}_{\beta}(S)$, we have the inequalities
\begin{displaymath} 2^{\alpha}4^{\beta - 1}\frac{|S|}{w(S)}
\leq {\Big[}\frac{\sigma(Q)}{|Q|}{\Big]^{p-1}} \leq 2^{\alpha + 1}4^{\beta}\frac{|S|}{w(S)}, \qquad Q \in \mathcal{P}_{\beta}(S). 
\end{displaymath}
Thus the measure $\sigma$ is not impossibly far away from Lebesgue measure on the cubes $Q \in \mathcal{P}_{\beta}(S)$, and this allows us to make use of the previously established estimate for the Lebesgue measure. As was already noted during the proof of \eqref{e.Lexp1}, any set of the form $\{x \in S : \phi_{S,\beta}^{\ast}(x) > \lambda\}$ can be expressed as the disjoint union of the \emph{maximal} cubes $Q \in \mathcal{D}_{\kappa}$ for which $Q \subset S$ and
\begin{equation}\label{de:form1} \left|\sum_{Q' \in \mathcal{P}_{\beta}(S) \colon Q \subsetneq Q' \subset S} \phi_{Q} \right| > \lambda. 
\end{equation}
Apply this with $\lambda = 20t \cdot w(S)/|S|$, and let $Q_{1},Q_{2},\ldots$ be the maximal cubes mentioned above.
{Note that it cannot happen that $Q_j=S$, since then the summation condition in \eqref{de:form1} would be empty ($S\subsetneq Q'\subset S$), and the estimate \eqref{de:form1} could not possibly hold. Thus} $Q_{j} \subsetneq S$ for every $j$, but there is no reason why $Q_{j} \in \mathcal{P}_{\beta}(S)$. Instead, $\tilde{Q}_{j} \in \mathcal{P}_{\beta}(S)$, where $\tilde{Q}_{j}$ denotes the parent (in $\mathcal{D}_{k})$ of $Q_{j}$. This follows from the maximality of $Q_{j}$ and the fact that the summation in \eqref{de:form1} is only over $Q' \in \mathcal{P}_{\beta}(S)$. Now let
\begin{displaymath} E_{\beta}(t) := \bigcup_{j = 1}^{\infty} \tilde{Q}_{j}. 
\end{displaymath}
This union may be assumed disjoint, since it is anyway the disjoint union of the maximal cubes among the cubes $\tilde{Q}_{j} \in \mathcal{P}_{\beta}(S)$. Recalling the estimate $|\phi_{Q}(x)| \leq 4^{-\beta + 1}w(S)|S|^{-1}$ used already in the proof of \eqref{e.Lexp1}, we may now estimate
\begin{align*} \left|\sum_{Q' \in \mathcal{P}_{\beta}(S) \colon \tilde{Q}_{j} \subsetneq Q' \subset S} \phi_{Q'}(x) \right| & \geq \left|\sum_{Q' \in \mathcal{P}_{\beta}(S) \colon Q_{j} \subsetneq Q' \subset S} \phi_{Q'}(x) \right| - |\phi_{\tilde{Q}_{j}}(x)|\\
& > 20t\frac{w(S)}{|S|} - 4 \cdot 4^{-\beta}\frac{w(S)}{|S|} \stackrel{(\ast)}{\geq} 16t\frac{w(S)}{|S|}, \qquad x \in \tilde{Q}_{j}. 
\end{align*} 
Of course, $(\ast)$ only holds in case $t \geq 4^{-\beta}$. One should also observe that we needed here the fact that the sum in the first expression is constant on $\tilde{Q}_{j}$, the $\mathcal{D}_{\kappa}$ child of the first cubes $Q'$ in the actual summation: this becomes essential, if $x \in \tilde{Q}_{j}$ is not in $Q_{j}$ to begin with, since we only have the inequality \eqref{de:form1} for $x \in Q_{j}$. The previous estimate now shows that $E_{\beta}(t) \subset \{x : \phi_{S,\beta}(x) > 16t\cdot w(S)/|S|\}$ for $t \geq 4^{-\beta}$, whence
\begin{displaymath} |E_{\beta}(t)| \leq \sqrt{2} \cdot 2^{-t4^{\beta}/\mathsf{K}}|S| \leq 2\cdot 2^{-t4^{\beta}/\mathsf{K}}|S| 
\end{displaymath}
by \eqref{de:form2}. For $t < 4^{-\beta}$ this estimate says nothing, so the same bound holds for all $t > 0$. Now sum over the disjoint cubes that form $E_{\beta}(t)$ and use $\sigma(Q) \leq {\big(} 2^{\alpha + 1}4^{\beta}w(S)^{-1}|S| {\big)^{1/(p-1)}}|Q|$ to obtain
\begin{equation*}
\begin{split}
   \sigma(E_{\beta}(t))
   &\leq {\Big(}2^{\alpha + 1}4^{\beta}\frac{|S|}{w(S)} {\Big)^{1/(p-1)}} |E_{\beta}(t)|
      \leq {\Big(}2^{\alpha + 1}4^{\beta}\frac{|S|}{w(S)} {\Big)^{1/(p-1)}}2 \cdot 2^{-t4^{\beta}/\mathsf{K}}|S| \\
   &\leq { (2\cdot 4^{\beta})^{1/(p-1)}\frac{\sigma(S)}{|S|}\cdot2 \cdot 2^{-t4^{\beta}/\mathsf{K}}|S|
     =2^{p'}\cdot 4^{\beta/(p-1)}}\cdot  2^{-t4^{\beta}/\mathsf{K}}\sigma(S).
\end{split}
\end{equation*}
In the last inequality we used again the fact that the cubes in $\mathcal{P}(S)$ are $(\sigma,\alpha,\kappa)$-adapted. The estimate for $\sigma(\{x \in S : \phi_{S}^{\ast}(x) > 20t \cdot w(S)/|S|\})$ is now obtained in the same manner as at the end of the proof of \eqref{e.Lexp1}:
\begin{align}
  \sigma(S)^{-1} & \left|\left\{x \in S : \phi_{S}^{\ast}(x) > 20t\frac{w(S)}{|S|}\right\}\right|\notag\\
  & \leq \sigma(S)^{-1}\sum_{\beta = 0}^{\infty} \left|\left\{x \in S : \phi_{S,\beta}^{\ast}(x) > 20 \cdot 2^{-\beta - 1}t\frac{w(S)}{|S|}\right\}\right|\notag\\
  &\label{form3} \leq {2^{p'}}\sum_{\beta = 0}^{\infty} {4^{\beta/(p-1)}}2^{-t2^{\beta - 1}/\mathsf{K}}
  = {2^{p'}}\sum_{\beta = 0}^{\infty} 2^{2\beta{/(p-1)} - t2^{\beta}/2\mathsf{K}}. 
\end{align} 
If $t < 2\mathsf{K}$, we clearly have
\begin{displaymath} \left|\left\{x \in S : \phi_{S}^{\ast}(x) > 20t\frac{w(S)}{|S|}\right\}\right| \leq 2 \cdot 2^{-t/2\mathsf{K}}\sigma(S). 
\end{displaymath}
We wish to prove a similar bound in case $t \geq 2\mathsf{K}$ by estimating the last series on line \eqref{form3}. Since
{ $2^{\beta} \geq \big(2/(p-1)+1\big)\beta + 1$ for $\beta \geq \beta_p$}, we have
\begin{equation*}
  2\beta/(p-1)-t2^{\beta}/2\mathsf{K}
  \leq -\beta-t/2\mathsf{K},\qquad\beta\geq\beta_p
\end{equation*}
Thus
\begin{equation*}
    2^{p'}\sum_{\beta = 0}^{\infty} 2^{2\beta/(p-1) - t2^{\beta}/2\mathsf{K}}
    \leq 2^{p'}\Big(\sum_{\beta=0}^{\beta_p-1} 2^{2\beta/(p-1)}+\sum_{\beta=\beta_p}^{\infty} 2^{-\beta}\Big)2^{-t/2\mathsf{K}}
    \leq C_p 2^{-t/2\mathsf{K}}.
\end{equation*}
Hence
\begin{displaymath} \left|\left\{x \in S : \phi_{S}^{\ast}(x) > 20t\frac{w(S)}{|S|}\right\}\right| \leq C_p \cdot 2^{-t/2\mathsf{K}}\sigma(S). 
\end{displaymath}
in each case. This completes the proof of \eqref{e.Lexp2}, and thus of Lemma~\ref{l.distribL}.
\section{Verifying the Testing Conditions} \label{s.testing}

The main Theorems of \S\ref{s.2wt} state that we can reduce the verification of the estimates of the main technical Theorem~\ref{t.haarShiftWtd} to this Lemma.  

\begin{lemma}\label{l.testing}  For an $ L ^2 $-bounded Haar shift operator $ \mathbb S $ with complexity $ \kappa $, 
a choice of $ 1<p< \infty $, and $ w \in A_p$, we have these estimates uniform over selection of dyadic cube $ Q$ and 
choice of linearization $ \mathbb L $ of the maximal truncations $ \mathbb S _{\natural}$.  
\begin{align}\label{e.testing1}
\lVert \mathbf 1_{Q} \mathbb L ^{\ast} (w \varphi  \mathbf 1_{Q})\rVert_{ L ^{p'} (\sigma )} 
& \lesssim   \kappa \lVert w\rVert_{A_p} w (Q) ^{1/p'} \,, 
\\ \label{e.testing3}
\Bigl\lVert  \sup_{R\subset Q}%
\mathbf{1}_{R} \frac{1}{w (R)}\int_{R}\left\vert \mathbb{L}^{\ast
}(\mathbf{1}_{R}\varphi w )(y)\right\vert \sigma (dy)  \Bigr\rVert_{ L ^{p}(w)} 
&\lesssim \kappa  \lVert w\rVert_{A_p} ^{1/(p-1)} \sigma(Q) ^{1/p} \,. 
\end{align}
In these estimates, $ p'= p / (p-1)$ is the conjugate index, $ \sigma = w ^{1-p'}$, and $ \varphi $ is a measurable 
function with $ \lVert \varphi \rVert_{\infty } \le 1$.   

\end{lemma}

We shall return to the proof of the Lemma above, turning here to the proofs of the two main technical estimates in Theorem~\ref{t.haarShiftWtd}.

\begin{proof}[Proof of the weak-type estimate \eqref{e.Sweak}.]  
We apply the inequality \eqref{e.TW2} from Theorem~\ref{weakFlat} giving sufficient conditions for the weak-type inequality  in the general two-weight case, and then Buckley's bound for the weak-type maximal inequality constant $\mathfrak{M}_{p,\operatorname{weak}}$ and the bound \eqref{e.testing1} for the `backward' testing constant $\mathfrak{T}_p$ in the one-weight case:
\begin{equation*}
\begin{split}
  \left\Vert \mathbb{S}_{\natural }(f\sigma )\right\Vert _{L^{p,\infty}(w )}
  &\lesssim \left( \kappa\mathfrak{M}_{p, \textup{weak}}+\mathfrak{T}_p\right) \left\Vert f\right\Vert _{L^{p}(\sigma )}\\
  &\lesssim\left(\kappa\|w\|_{A_p}^{1/p}+\kappa\|w\|_{A_p}\right)\|f\|_{L^p(\sigma)}
  \lesssim\kappa\|w\|_{A_p}\|f\|_{L^p(\sigma)},
\end{split}
\end{equation*}
noting that the second term dominates since $\|w\|_{A_p}\geq 1$.
\end{proof}

\begin{proof}[Proof of the strong-type estimate \eqref{e.Sstrong}.] 
We apply the sufficient conditions in the 
two weight setting of Theorem~\ref{t.genlStrong}, and then the available estimates for the testing constants in the one-weight situation: Buckley's bound for the strong-type maximal inequality constant $\mathfrak{M}_p$, and the bounds \eqref{e.testing1} and \eqref{e.testing3} for the `backward' and nonstandard testing constants $\mathfrak{T}_p$ and~$\mathfrak{N}_p$:
\begin{equation}
\begin{split}
  \lVert \mathbb{S}_{\natural }(f\sigma )\rVert _{L^{p}(w )}
  &\lesssim\left\{ \kappa \mathfrak{M}_{p}
     +\mathfrak{T}_p
      +\mathfrak{N}_{p}\right\} \lVert f\rVert _{L^{p}(\sigma )}\\
   &\lesssim\left\{\kappa\|w\|_{A_p}^{1/(p-1)}
     +\kappa\|w\|_{A_p}+\kappa\|w\|_{A_p}^{1/(p-1)}\right\}\|f\|_{L^p(\sigma)} \\
    &\lesssim\kappa\|w\|_{A_p}^{\max\{1,1/(p-1)\}}\|f\|_{L^p(\sigma)}.
\end{split}
\end{equation}%
In this case either the first and third, or the second term dominates, depending on whether $p<2$ or $p>2$.
%
%
\end{proof}

We turn to the proof of lemma~\ref{l.testing}.

\begin{proof}[Proof of \eqref{e.testing1}.] 
The data is fixed: of index $ p$, weight $ w \in A_p$,  cube $Q $, function $\varphi  $  bounded by $ 1$, and Haar Shift $ \mathbb S $ of complexity $ \kappa $.  

The definition of $ \mathbb L ^{\ast} $ is as in \eqref{e.adjoint}, and is in particular a sum over all dyadic cubes $ P$.  Now, the sum in \eqref{e.adjoint} can obviously be restricted to a sum over $ P$ that intersect $ Q$.  Moreover, the sum over $ P$ that contain $ Q$ can be controlled, using a straight forward appeal to the size conditions on the Haar Shift operators, and the definition of $ A_p$.

Thus, in what follows, we need only consider cubes which are contained in $ Q$.
We split these cubes into $ (\sigma ,\alpha , \kappa )$-adapted subcollections, in the sense of Definition~\ref{d.adapted}. Namely, 
 we form the $\kappa+1$ subcollections according to the value of $\log_2\ell(Q)\mod(\kappa+1)$. And each of them is further split according to the unique number $\alpha$ such that
\begin{equation*}
  2^{\alpha}\leq\frac{w(Q)}{|Q|}\Big(\frac{\sigma (Q)}{|Q|}\Big)^{p-1}<2^{\alpha+1},
\end{equation*}
 where $0\le \alpha \le \lceil\log_2 \lVert w \rVert_{A_{p}}\rceil$ is an integer.   

For any one of these $(\sigma,\alpha,\kappa)$-adapted subcollections, say $\mathcal{Q}$, apply the Definition~\ref{d.stopping}, to get $ w$-stopping cubes $ \mathcal S$, and corona decomposition $ \{\mathcal P (S) \;:\; S\in \mathcal S\}$.   Define functions 
$
\mathbb F  _{S} := \mathbb L ^{\ast} _{\mathcal P (S)} (w\varphi \mathbf 1_{S}) 
$.
For these functions, we have the distributional estimate \eqref{e.Lexp2}, and it remains to show that 
\begin{equation}\label{e.F1<}
\Bigl\lVert  \sum_{S\in \mathcal S} \mathbb F_{S}  \Bigr\rVert_{ L ^{p'} (\sigma )} \lesssim 
  2 ^{\alpha/{p}} \lVert w\rVert_{A_p} ^{1/p'} w (Q) ^{1/p'} .
\end{equation}
To conclude \eqref{e.testing1}, we sum this estimate over $\alpha$ to get the upper bound
$
\|w\|_{A_p} w(Q)^{1/p'},
$, 
and then over the $\kappa+1$ possible values of $\log_2\ell(Q')\mod(\kappa+1)$ to get the estimate \eqref{e.testing1} as stated, with the factor $\kappa$.

The distributional estimates are exponential in nature, which permit a facile estimate of this norm.  Define level sets of the functions $ \mathbb F _S$ by 
\begin{gather}\label{e.F10}
F _{S,0} :=\bigl\{\lvert  F_S\rvert < \mathsf K  \tfrac {w (S)} {\lvert  S\rvert}  \bigr\}\,, 
\\
\label{e.F1n}
F _{S,n} :=\bigl\{ \mathsf K n   \tfrac {w (S)} {\lvert  S\rvert} \le    \lvert  F_S\rvert < \mathsf K (n+1)   \tfrac {w (S)} {\lvert  S\rvert}  \bigr\}  \,, 
\qquad n\ge 1 \,,
\end{gather}
where $\mathsf K$ is the dimensional constant in \eqref{e.Lexp2}.
 And then we estimate, by a familiar trick, 
\begin{align}
\Bigl\lvert  \sum_{S\in \mathcal S} \mathbb F  _{S}  \Bigr\rvert ^{p'} & 
= 
\Bigl\lvert  \sum_{n=0} ^{\infty } \sum_{S\in \mathcal S}   (1+n) ^{-2/p + 2/p} \mathbb F  _{S} \mathbf 1_{F_ {S,n}}    \Bigr\rvert ^{p'} 
\\
& 
\le \big(\frac {\pi ^2 } 2\big)^{p'/p}   \sum_{n=0} ^{\infty }  (1+n) ^{2p'/p} 
\Bigl\lvert \sum_{S\in \mathcal S}    \mathbb F  _{S}   \mathbf 1_{F_ {S,n}}  \Bigr\rvert ^{p'} \\
& \le \big(\frac {\pi ^2 } 2\big)^{p'/p}  \mathsf K^{p'} \sum_{n=0} ^{\infty }  (1+n) ^{2p'/p+p'}\Big(\sum_{S\in\mathcal{S}}\frac{w(S)}{|S|}\mathbf 1_{F_{S,n}}\Big)^{p'}
\\ \label{e.familiar}
&\le  \big(\frac {\pi ^2 } 2\big)^{p'/p}   \mathsf K ^{p'} \sum_{n=0} ^{\infty }  (1+n) ^{2p'/p +p'} 
  {\Big(\frac{4}{3}\Big)^{p'}}\sum_{S\in \mathcal S} \Bigl[\frac {w (S)} {\lvert  S\rvert } \Bigr] ^{p'} 
 \mathbf 1_{F_ {S,n}}.
\end{align}
Here, we use H\"older's inequality in the first inequality, and in the fourth, we use the fact that the averages 
$   \tfrac {w (S)} {\lvert  S\rvert} $ increase geometrically along the stopping intervals.   This is wasteful in the powers on $ n$, 
which is of no consequence for us.  

We then prove \eqref{e.F1<} as follows, absorbing the dimensional constant $\mathsf K$ into the constants implicit in the notation $\lesssim$. We have 
\begin{align}
\Bigl\lVert  \sum_{S\in \mathcal S} \mathbb F_{S}  \Bigr\rVert_{ L ^{p'} (\sigma )} ^{p'} 
& 
\lesssim 
 \sum_{n=0} ^{\infty }  (1+n) ^{2p'/p +p'}  \sum_{S\in \mathcal S}   \Bigl[\frac {w (S)} {\lvert  S\rvert } \Bigr] ^{p'}  
 \sigma (F _{S,n}) 
 \\  
 & \lesssim  
  \sum_{n=0} ^{\infty }  (1+n) ^{2p'/p +p'}   2 ^{-n} \sum_{S\in \mathcal S}   \Bigl[\frac {w (S)} {\lvert  S\rvert } \Bigr] ^{p'}  \sigma (S) 
  \\ \label{e.veryS}
 & \lesssim  
    \sum_{S\in \mathcal S}  w (S)  \Bigl[\frac {w (S)} {\lvert  S\rvert } \Bigr] ^{p' -1} \frac  {\sigma (S)} {\lvert  S\rvert }   \\
 & \lesssim   \label{eq:ApStopping} 
   2 ^{ (p'-1)\alpha}\sum_{S\in \mathcal S} w (S) 
  \lesssim  
   2 ^{ (p'-1)\alpha}  \lVert w\rVert_{A_p} w (Q) \,. 
\end{align}
To pass to the second line, we use the distributional estimate \eqref{e.Lexp2}, giving us a sum in $ n$ that is trivially bounded. 
From the third to the fourth line, we use the condition \eqref{e.ApFixed} of Definition~\ref{d.adapted} and then use the estimate \eqref{e.stoppingSum}. 
\end{proof}

\begin{proof}[Proof of \eqref{e.testing3}]
Fix the relevant data for this condition, and take a dyadic cube $ R  \subset Q$.  
We take $ \widetilde{\mathbb L} ^{\ast } _{R} $ to be as in \eqref{e.adjoint}, but with the sum restricted to cubes that 
contain $ R$.   Then, 
\begin{equation*}
\sup _{R\subset Q} \frac 1 {w (R)} \int _{R} \lvert   \widetilde{\mathbb L} ^{\ast } _{R} (w \varphi \mathbf 1_{R})\rvert  \; \sigma (dy) 
\le\sup _{R\subset Q} \frac {\sigma (R)} {\lvert  R\rvert } \le  M (\sigma \mathbf 1_{Q}) \,.  
\end{equation*}
And this last term is controlled by Buckley's estimate for the Maximal Function.  
Thus, when we consider the expression 
\begin{equation*}
 \frac{1}{w (R)}\int_{R}\left\vert \mathbb{L}^{\ast
}(\mathbf{1}_{R}\varphi w )(y)\right\vert \sigma (dy) 
\end{equation*}
we can in addition assume that all cubes contributing to $ \mathbb L ^{\ast}  $ are contained in $ R$, namely we can replace 
$ \mathbb L ^{\ast}  $ by $ \mathbb L ^{\ast} _{\mathcal R (R)}$, where $ \mathcal R (R)$ is an appropriate collection of cubes contained in $ R$. 

To make a relevant estimate of this integral, we can consider $ \mathcal Q$ a collection of cubes $ Q' \subset Q$ that are 
$ (w,\alpha , \kappa )$-adapted, where  $ 0\le \alpha \le \lceil\log_2 \lVert w \rVert_{A_{p}}\rceil$.  Then, we will show that 
\begin{equation}\label{e.T3<}
  \Bigl\lVert  \sup _{R \subset Q}  \mathbf 1_{R}
    \frac{1}{w (R)}\int_{R}\left\vert \mathbb{L}^{\ast} _{\mathcal R (R)}(\mathbf{1}_{R}\varphi w )(y)\right\vert \sigma (dy) \Bigr\rVert_{ L ^{p} (w)} 
    \lesssim 
    2 ^{\alpha / (p-1)}  \sigma (Q) ^{1/p} \,, 
\end{equation}  
where 
the collections $ \mathcal R (R)$ are those cubes $ Q' \in \mathcal Q$ 
with $ Q'\subset R$.   Summing this estimate over  $ 0\le \alpha \le \lceil\log_2 \lVert w \rVert_{A_{p}}\rceil$ and the $\kappa+1$ possible values of $\log_2\ell(Q')\mod(\kappa+1)$ will prove \eqref{e.testing3}. 

We then apply Definition~\ref{d.stopping}, letting $ \mathcal  S (R)$ be a collection of stopping cubes for $ \mathcal R (R)$, with 
Corona Decomposition $ \{\mathcal P (S) \;:\; S\in \mathcal S (R)\}$.  Then, we have 
\begin{align}
 \frac{1}{w (R)}\int_{R}\left\vert \mathbb{L}^{\ast
}_{\mathcal R (R)}(\mathbf{1}_{R}\varphi w )(y)\right\vert \sigma (dy)  
& \le
 \frac{1}{w (R)}\sum_{S\in \mathcal S (R)} 
\int_{R}\left\vert \mathbb{L}^{\ast
} _{\mathcal P (S)}(\mathbf{1}_{R}\varphi w )(y)\right\vert \sigma (dy)  
\\
& \lesssim  \mathsf K \frac{1}{w (R)} \sum _{S\in \mathcal S (R)} \frac {w (S)} {\lvert  S\rvert } \sigma (S) 
\\
 &  \lesssim \mathsf K  2 ^{\alpha/ (p-1) } \frac{1}{w (R)}
    \sum _{S\in \mathcal S (R)} w(S) \Bigl[ \frac {\lvert  S\rvert } {w (S)} \Bigr] ^{1/ (p-1)} \,. 
\end{align}
Here, we use the estimate \eqref{e.Lexp2}. Recalling that $\mathsf K$ is only a dimensional constant, we henceforth absorb it from the notation.  And so, to prove \eqref{e.T3<}, we should show that 
\begin{equation}\label{e.T3<<}
  \Bigl\lVert  
\sup _{R \subset Q} \mathbf 1_{R} \frac{1}{w(R)}
  \sum _{S\in \mathcal S (R)} w(S) \Bigl[ \frac {\lvert  S\rvert } {w (S)} \Bigr] ^{1/ (p-1)} 
\Bigr\rVert_{ L ^{p} (w)} \lesssim  \sigma (Q) ^{1/p} \,.
\end{equation}

We estimate
\begin{equation*}
\begin{split}
  \sum _{S\in \mathcal S (R)}  w(S)\Bigl[ \frac {\lvert  S\rvert } {w (S)} \Bigr] ^{1/ (p-1)} 
  & =\sum_{k=0}^{\infty}\sum _{S\in \mathcal S_k (R)}  w(S)\Bigl[ \frac {\lvert  S\rvert } {w (S)} \Bigr] ^{1/ (p-1)}  \\
  & \leq\sum_{k=0}^{\infty}\sum _{S\in \mathcal S_k (R)}  w(S)\Bigl[4^{-k} \frac {\lvert  R\rvert } {w (R)} \Bigr] ^{1/ (p-1)}  \\ 
  & \leq\sum_{k=0}^{\infty}w(R)\Bigl[4^{-k} \frac {\lvert  R\rvert } {w (R)} \Bigr] ^{1/ (p-1)}\\
 &  \lesssim w(R)\Bigl[ \frac {\lvert  R\rvert } {w (R)} \Bigr] ^{1/ (p-1)}  = w(R)\bigl[\mathbb E _{R} ^{w} (w ^{-1}\mathbf 1_Q ) \bigr] ^{1/ (p-1)} \,, 
\end{split}
\end{equation*}
where we used the fact that the terms 
$  \frac {\lvert  S\rvert } {w (S)}$ decrease geometrically along the levels of the stopping cubes.  Therefore, 
 to see the estimate \eqref{e.T3<<}, we should estimate 
 \begin{displaymath}
 \lVert  M _{w} (w ^{-1} \mathbf 1_{Q} ) ^{ 1/ (p-1)} \rVert_{ L ^{p} (w)} 
 = \lVert  M _{w} (w ^{-1} \mathbf 1_{Q} ) \rVert_{ L ^{p'} (w)}  ^{ 1/ (p-1)}
 \lesssim \lVert w ^{-1} \mathbf 1_{Q} \rVert_{ L ^{p'} (w)}  ^{ 1/ (p-1)} 
 = \sigma (Q) ^{1/p}
 \end{displaymath}
by the universal maximal function estimate.  
\end{proof}

\section{Variations on the Main Theorem}

In this final section, we present a couple of variations of the main result, Theorem~\ref{t.main}, where we allow the appearance of different weight characteristics than just $\|w\|_{A_p}$ in the norm estimates in $L^p(w)$. First, as a direct consequence of Theorem~\ref{t.main}, we deduce the following result, which was conjectured by Lerner and Ombrosi \cite[Conjecture 1.3]{LO} for untruncated operators $T$:

\begin{corollary}
For $  T $ an $ L ^{2} (\mathbb R ^{d})$ bounded  Calder\'on-Zygmund Operator and $1<q<p<\infty$, 
\begin{equation*}
  \|T_{\natural}f\|_{L^p(w)}\leq C_{T,p,q}\|w\|_{A_q}\|f\|_{L^p(w)}.
\end{equation*}
\end{corollary}

\begin{proof}
Choose some $r\in(p,\infty)$. Since $\|w\|_{A_r}\leq\|w\|_{A_q}$, the weak-type estimate of Theorem~\ref{t.main} shows that
\begin{equation*}
\begin{split}
  \|T_{\natural}f\|_{L^{q,\infty}(w)}
  &\leq C_{T,q}\|w\|_{A_q}\|f\|_{L^q(w)}, \\
  \|T_{\natural}f\|_{L^{r,\infty}(w)}
  &\leq C_{T,r}\|w\|_{A_r}\|f\|_{L^r(w)}
    \leq C_{T,r}\|w\|_{A_q}\|f\|_{L^r(w)}.
\end{split}
\end{equation*}
It suffices to apply the Marcinkiewicz interpolation theorem to the sublinear operator $T_{\natural}$ to deduce the asserted strong-type bound in $L^p(w)$. (For $p\geq 2$, we could also have used the strong-type estimate of Theorem~\ref{t.main} and $\|w\|_{A_p}\leq\|w\|_{A_q}$, without any interpolation.)
\end{proof}

In order to describe the other variant of the main theorem, we consider the \emph{$A_\infty$ characteristic}
\begin{equation}\label{eq:AinftyNorm}
  \|w\|_{A_\infty}:=\sup_Q\frac{1}{w(Q)}\int_Q M(w\mathbf{1}_Q),
\end{equation}
which has been introduced (with a different notation) by Wilson \cite{MR883661}, and more recently used by Lerner \cite{1005.1422}. One can show (see \cite{HytPer}) that $\|w\|_{A_\infty}\leq c_d\|w\|_{A_p}$ for all $p\in(1,\infty)$, and therefore the following estimate is an improvement of Theorem~\ref{t.main}:

\begin{theorem}\label{t.CZopAinfty}
For $  T $ an $ L ^{2} (\mathbb R ^{d})$ bounded  Calder\'on-Zygmund Operator and $1<p<\infty$, 
\begin{align}
\lVert   T _{\natural } f  \rVert _{ L ^{p,\infty } (w)} &\le C_{T,p} \|  w \|_{ A _{p}}^{1/p} \|  w \|_{ A _{\infty}}^{1/p'} \lVert  f \rVert _{L ^{p} (w)},\\
\lVert   T _{\natural} f  \rVert _{ L ^{p } (w)} &\le C_{T,p}\big( \|  w \|_{ A _{p}}^{1/p} \|  w \|_{ A _{\infty}}^{1/p'}+\|w\|_{A_p}^{1/(p-1)}\big) \lVert  f \rVert _{L ^{p} (w)}.
\end{align}
\end{theorem}

Note that the strong-type bound is a strict improvement of the corresponding bound in Theorem~\ref{t.main} only for $2<p<\infty$. For $1<p\leq 2$, the term $\|w\|_{A_p}^{1/(p-1)}$ dominates, and this is exactly the same bound as in Theorem~\ref{t.main}.

An earlier result of this type is implicitly contained in the work of Lerner \cite[Section 5.5]{1005.1422}:
\begin{equation*}
  \|T_{\natural}f\|_{L^p(w)}\leq C_{T,p}\|w\|_{A_p}^{1/2}\|w\|_{A_\infty}^{1/2}\|f\|_{L^p(w)},\qquad 3\leq p<\infty.
\end{equation*}
Our estimate improves on this, as both terms $\|w\|_{A_\infty}^{1/(p-1)}\leq\|w\|_{A_\infty}^{1/2}$ (for $p-1\geq 2$) and
\begin{equation*}
   \|  w \|_{ A _{p}}^{1/p} \|  w \|_{ A _{\infty}}^{1/p'}
   =\|  w \|_{ A _{p}}^{1/2} \|  w \|_{ A _{\infty}}^{1/2}\Big(\frac{\|w\|_{A_\infty}}{\|w\|_{A_p}}\Big)^{1/2-1/p}
\end{equation*}
are smaller than $\|w\|_{A_p}^{1/2}\|w\|_{A_\infty}^{1/2}$ for $3\leq p<\infty$.
 
On the other hand, Theorem~\ref{t.CZopAinfty} fails to reproduce the following bound for untruncated operators $T$, which was recently obtained in \cite{HytPer}:
\begin{equation*}
  \|Tf\|_{L^2(w)}\lesssim \|w\|_{A_2}^{1/2}\big(\|w\|_{A_\infty}+\|w^{-1}\|_{A_\infty}\big)^{1/2}\|f\|_{L^2(w)}.
\end{equation*}
We do not know whether this bound remains true for $T_{\natural}$ in place of $T$. More generally, although the optimal $L^p(w)$ bounds in terms of $\|w\|_{A_p}$ coincide for $T$ and $T_{\natural}$, we do not know if this is the case when allowing the more complicated dependence on both $A_p$ and $A_\infty$ characteristics.

We now discuss the proof of Theorem~\ref{t.CZopAinfty}. 
Just like the proof of Theorem~\ref{t.main}, thanks to the Representation Theorem~\ref{t.bcr}, this is a consequence of an analogous result for Haar shifts:

\begin{theorem}\label{t.shiftAinfty}
Let $ \mathbb S $ be a Haar shift operator with complexity $ \kappa $, a paraproduct, or a dual paraproduct.
For $ 1< p < \infty $ and $ w \in A_p$,  we have the estimates 
\begin{align}
\lVert \mathbb S_{\natural } f \rVert_{ L ^{p, \infty } (w)} &\lesssim \kappa \|  w \|_{ A _{p}}^{1/p} \|  w \|_{ A _{\infty}}^{1/p'} \lVert f\rVert_{L ^{p} (w)}  \,, 
\\ 
\lVert \mathbb S_{\natural } f \rVert_{ L ^{p } (w)} & \lesssim \kappa \big( \|  w \|_{ A _{p}}^{1/p} \|  w \|_{ A _{\infty}}^{1/p'}+\|w\|_{A_p}^{1/(p-1)}\big)
   \lVert f\rVert_{L ^{p} (w)}.
\end{align}
\end{theorem}

\begin{proof}[Sketch of proof]
We can exploit the same two-weight estimates as before,
\begin{align}
  \lVert \mathbb S_{\natural }(f\sigma) \rVert_{ L ^{p, \infty } (w)}
  &\lesssim \left\{\kappa\mathfrak{M}_{p,\operatorname{weak}}+\mathfrak{T}_p\right\}\lVert f\rVert_{L ^{p} (\sigma)}, \\
  \lVert \mathbb S_{\natural } (f\sigma) \rVert_{ L ^{p } (w)}
  &\lesssim \left\{\kappa\mathfrak{M}_{p}+\mathfrak{T}_p+\mathfrak{N}_p\right\}\lVert f\rVert_{L ^{p} (\sigma)},
\end{align}
and it only remains to see how the $A_\infty$ norm can be incorporated into the estimates of the testing constants appearing here. Recall that we have
\begin{equation*}
  \mathfrak{M}_{p,\operatorname{weak}}
  \lesssim\|w\|_{A_p}^{1/p},\qquad
  \mathfrak{M}_{p}\lesssim\|w\|_{A_p}^{1/(p-1)},\qquad
  \mathfrak{N}_p\lesssim\kappa\|w\|_{A_p}^{1/(p-1)},
\end{equation*}
which are of admissible size for Theorem~\ref{t.shiftAinfty}. The improvement over the earlier bounds comes from a more careful estimation of $\mathfrak{T}_p$.

Recall that the estimate $\mathfrak{T}_p\lesssim\kappa\|w\|_{A_p}$ is proven in Lemma~\ref{l.testing}, Eq. \eqref{e.testing1}, and that this proof is reduced to the estimate \eqref{e.F1<}. The final step in the proof of this estimate, \eqref{eq:ApStopping}, is an application of Lemma~\ref{l.stoppingSum}, which says that the collection $\mathcal{S}$ of $w$-stopping cubes for a cube $Q$ satisfies
\begin{equation*}
  \sum_{S\in\mathcal{S}}w(S)\lesssim\|w\|_{A_p}w(Q).
\end{equation*}
It is this last bound where $\|w\|_{A_p}$ can be replaced by $\|w\|_{A_\infty}$, as observed in \cite{HytPer}:
\begin{equation*}
  \sum_{S\in\mathcal{S}}w(S)\lesssim\|w\|_{A_\infty}w(Q).
\end{equation*}
Using this bound in \eqref{eq:ApStopping}, we obtain
\begin{equation*}
  \mathfrak{T}_p\lesssim\kappa\|w\|_{A_p}^{1/p}\|w\|_{A_\infty}^{1/p'},
\end{equation*}
and thus the estimates as asserted in Theorem~\ref{t.shiftAinfty}.
\end{proof}

\begin{remark}
It is shown in \cite{HytPer} that there is the sharper bound $\mathfrak{M}_p\lesssim\big(\|w\|_{A_p}\|\sigma\|_{A_\infty}\big)^{1/p}$, where $\sigma=w^{-1/(p-1)}$. However, this does not allow us improve on the above estimates, as our bound for $\mathfrak{N}_p$, namely $\|w\|_{A_p}^{1/(p-1)}$, is already as large as Buckley's classical bound for $\mathfrak{M}_p$.
\end{remark}

\begin{remark}
Another, perhaps more commonly used definition of the $A_\infty$ characteristic is the following quantity introduced by Hru\v{s}\v{c}ev \cite{MR727244}:
\begin{equation*}
  \sup_Q\Big(\frac{1}{|Q|}\int_Q w\Big)\exp\Big(\frac{1}{|Q|}\int_Q\log w^{-1}\Big).
\end{equation*}
The $A_\infty$ characteristic as defined in \eqref{eq:AinftyNorm} is always smaller (up to a dimensional constant), and it can be much smaller for some weights (see \cite{HytPer} for details), so Theorems~\ref{t.CZopAinfty} and \ref{t.shiftAinfty} are sharper when stated in terms of $\|w\|_{A_\infty}$ as in \eqref{eq:AinftyNorm}.
\end{remark}

\end{document}